\documentclass[notitlepage,twoside,a4paper,8.5pt]{amsart}
\usepackage{setspace}
\usepackage{graphicx}
\usepackage[dvipsnames]{xcolor}
\usepackage{mathrsfs}
\usepackage{amsmath,amssymb,enumerate}
\usepackage{epsfig,fancyhdr,color}
\usepackage{epstopdf}
\usepackage{amssymb}
\usepackage{amsmath,amsthm}
\usepackage{latexsym}
\usepackage{amscd}
\usepackage{psfrag}
\usepackage{url}
\usepackage{stfloats}

\usepackage{epsf}
\usepackage[latin1]{inputenc}
\usepackage[all]{xy}
\usepackage{tikz}
\usepackage{prettyref}
\usepackage{subcaption}
\usepackage{float}
\usepackage{color}
\usepackage{array}
\usepackage{cite}
\usepackage[totalwidth=16cm,totalheight=23cm]{geometry}

\newcommand{\II}{I\hspace{-0.1cm}I}
\newcommand{\III}{I\hspace{-0.1cm}I\hspace{-0.1cm}I}

\makeatletter
\newcommand{\myitem}[1]{%
	\item[#1]\protected@edef\@currentlabel{#1}%
}
\makeatother

\newcommand{\f}{\mathsf{F}}
\newcommand{\g}{\mathsf{G}}
\newcommand{\fp}{\mathsf{F_{+}}}

\newcommand{\fm}{\mathsf{F_{-}}}

\newcommand{\fpm}{\mathsf{F_{\pm}}}

\newcommand{\ddt}{\frac{d}{dt}\Bigr|_{\substack{t=0}}}
\newcommand{\ddtt}{\frac{d}{ds}\Bigr|_{\substack{s=0}}}

\newcommand{\red}[1]{\textcolor{black}{#1}}
\newcommand{\deter}{\mathsf{det}}
\newcommand{\T}{\mathcal{T}(S)}
\newcommand{\QF}{\mathcal{QF}(S)}
\newcommand{\MF}{\mathcal{MF}(S)}
\newcommand{\FMF}{\mathcal{FMF}(S)}
\newcommand{\AF}{\mathcal{AF}(S)}
\newcommand{\F}{\mathcal{F}(S)}
\newcommand{\Wfp}{\mathcal{W}^{+}_{\mathsf{F}_{+}}}
\newcommand{\Wfm}{\mathcal{W}^{-}_{\mathsf{F}_{-}}}
\newcommand{\Wf}{\mathcal{W}_{\f}}
\newcommand{\Sig}{S}

\newcommand{\ML}{\mathcal{ML}(S)}

\newcommand{\p}{\mathsf{p}}
\newcommand{\pp}{\mathsf{P}}

\newcommand{\ext}{\mathsf{ext}}

\DeclareMathOperator{\trace}{tr}

\DeclareMathOperator{\hess}{Hess}

\newtheorem{theorem}{\rm\bf Theorem}[section]
\newtheorem{proposition}[theorem]{\rm\bf Proposition}
\newtheorem{lemma}[theorem]{\rm\bf Lemma}
\newtheorem{corollary}[theorem]{\rm\bf Corollary}
\newtheorem{definition}[theorem]{\rm\bf Definition}
\newtheorem{remark}[theorem]{\rm\bf Remark}

\def\interieur#1{\mathord{\mathop{\kern 0pt #1}\limits^\circ}}

\usepackage{hyperref}
\hypersetup{hidelinks}

\begin{document}

	\setcounter{tocdepth}{1}

	\title{Measured foliations at infinity of quasi-Fuchsian manifolds close to the Fuchsian locus }
	
\author[Diptaishik Choudhury]{Diptaishik Choudhury}
\address{Diptaishik Choudhury: Department of Mathematics, University of Luxembourg, 6 Avenue de la Fonte, Esch-sur-Alzette, Luxembourg- 4364} \email{diptaishik.choudhury@uni.lu}
\maketitle
	
	\begin{abstract}
		Measured foliations at infinity of quasi-Fuchsian manifolds are a natural analogue at infinity to the bending laminations on the boundary of its convex core. We show that given a pair of arational measured foliations $(\f_{+},\f_{-})$ which fill a closed hyperbolic surface $S$, for $t>0$ sufficiently small, $t\fp$ and $t\fm$ can be uniquely realised as the {measured foliations at infinity} of a quasi-Fuchsian manifold homeomorphic to $S\times \mathbb{R}$, which is sufficiently close to the Fuchsian locus. Here arationality means that the corresponding measured laminations are maximal. The proof is inspired by Bonahon's in\cite{bonahon05} which shows that a quasi-Fuchsian manifold close to the Fuchsian locus can be uniquely determined by the data of filling measured bending laminations on the boundary of its convex core. We also interpret the result in {half-pipe} geometry.
		
	\end{abstract}
\tableofcontents

	\section{Introduction}	
	\hspace*{.07 cm}  Let $S$ be a closed, oriented surface with genus $g \geq 2$ and $M$ a $3$-manifold homeomorphic to $S\times \mathbb{R}$. Call the space of isotopy classes of Fuchsian metrics on $M$ as the {Fuchsian locus} $\F$ and note that it can also be identified with the Teichm{\"u}ller space $\T$ (see \S \S  \ref{teichspace}). Let $T^{*}_{[c]}\T$ be its cotangent space at a point $[c]\in \T$ and again identify it with the space of holomorphic quadratic differentials $Q(S,[c])$ on $(S,[c])$ (see \S \S \ref{HQDintro}). Now, consider quasi-Fuchsian hyperbolic metrics on $M$ and let $\QF$ denote the space of isotopy classes of quasi-Fuchsian metrics on $M$. Denote the connected components of the boundary at infinity of $M$ as $\partial^{+}_{\infty}M$ and $\partial^{-}_{\infty}M$ (both being homeomorphic to $S$) and let $([c_{+}],[c_{-}])\in \mathcal{T}(\partial^{+}_{\infty}M)\times \mathcal{T}({\partial^{-}_{\infty}M})$ be the respective conformal classes (see Theorem \ref{bersthm}).
	\\ \hspace*{.5 cm} There are unique holomorphic maps, well-defined up to right composition by M{\"o}bius transformations, from the universal covers $\widetilde{\partial^{+}_{\infty}M}, \widetilde{\partial^{-}_{\infty}M}\subset \partial_{\infty}\mathbb{H}^{3}\cong \mathbb{C}P^{1}$ to the unit disc $\Delta\subset \mathbb{C}$ that we obtain by uniformising the respective complex structures (see \S \S \ref{introductionary}). Let the Schwarzians at infinity $\sigma_{+}\in Q(\partial^{+}_{\infty}M,[c_{+}])$ and $\sigma_{-}\in Q(\partial^{-}_{\infty}M,[c_{-}])$ be the holomorphic quadratic differentials obtained by taking the Schwarzian derivative of these maps respectively and passing to quotients. We define the measured foliations at infinity $\fp,\fm$ of a quasi-Fuchsian manifold $M$ as the horizontal measured foliations of $\sigma_{+},\sigma_{-}$ on $(\partial^{+}_{\infty}M,[c_{+}]),(\partial^{-}_{\infty}M,[c_{-}])$ respectively. These measured foliations at infinity can be seen as a natural analog at infinity to the bending lamination on the boundary of the convex core of a quasi-Fuchsian manifold (see \S\S \ref{schwatinf}, Lemma \ref{graft}). \\
	\hspace*{.5 cm} Let $\MF$ denote the space of equivalence classes of measured foliations on $S$  (see  \S \S\ref{defmf}, \S \S \ref{horifoli},\cite{AST_1979__66-67_}) and $\fp,\fm \in \MF$. Further, given a pair of measured foliations $(\f,\g)\in \MF\times \MF$ we have the notion of them being a pair which {fills} $\Sig$. That is to say, any other measured foliations $\mathsf{H}$ has non-zero intersection with both $\f$ or $\g$ (see Definition \red{\ref{fillupdef}} and \red{\S\S\ref{gardinermasur}}).
	So we ask (see Question $7.4$ in\cite{schlenker20}) whether is it possible to determine a quasi-Fuchsian manifolds uniquely by its measured foliations at infinity? \\
	\hspace*{.5 cm}Now let $\mathcal{MF}_{0}(S)\subset \MF$ be the subspace of measured foliations which are arational, i.e, all the prongs are of order $3$ and there are no leaves joining the prongs (see Definition \ref{arationaldef} and Lemma \ref{arationaliff}); $\FMF$ the space of all pairs of measured foliations that fill $S$ and $\mathcal{FMF}_{0}(S)$ be the subspace of such pair which are arational. If the pair $(\fp,\fm)$ belongs to $ \mathcal{FMF}_{0}(S)$, then so do the pair $(t\fp,t\fm)$, for all $t>0$ (see \S \S \ref{defmf}). Note also that for a metric $g\in \F$ the Schwarzians and the measured foliations at infinity  are zero (see \S \S \ref{introductionary}). The result of principal interest that answers the above question partially for quasi-Fuchsian manifolds near the Fuchsian locus is:
	{\theorem \label{thm1.1}
		For every pair of measured foliations $(\fp,\fm)$ which are arational and fill $S$, there exists an $\epsilon_{\fpm}>0$ such that for $\forall t \in (0,\epsilon_{\fpm})$ there exists an unique quasi-Fuchsian metric $g\in \QF$ on $M$ sufficiently close to the Fuchsian locus, whose measured foliations at infinity are given by  $t\fp$ and $t\fm$ .}\\\\
	That is, given the map $\mathfrak{F}: \QF \rightarrow \mathcal{MF}(\partial^{+}_{\infty}M) \times \mathcal{MF}({\partial^{-}_{\infty}M)}$ sending a quasi-Fuchsian metric to the measured foliations at infinity at the positive and negative end respectively; we have a unique solution $g\in \QF$ to the equation $\mathfrak{F}(g)=(t\fp,t\fm)$ when restricted to $(\fp,\fm)\in \mathcal{FMF}_{0}(S)$ for $t>0$ small enough. Now let $q^{\mathsf{H}}_{[c]}$ be the unique holomorphic quadratic differentials realising  $\mathsf{H}\in \MF$ as its horizontal measured foliation on $(S,[c])$ (see \S \S \ref{HMsection}). An immediate consequence along the lines of McMullen's quasi-Fuchsian reciprocity (see\cite{mcmullen,Krasnov2009,filippoksurface}) which helps in describing the Schwarzians at infinity is that if $g$ be a quasi-Fuchsian metric on $M$ such that the measured foliations at infinity are given as $(t\fp,t\fm)$ for some filling arational pair $(\fp,\fm)$, $t>0$ small enough; then the Schwarzians at infinity of $(M,g)$ are $t^{2}q^{\fp}_{[c_{+}]}\in Q(\partial^{+}_{\infty}M,[c_{+}])$ and $t^{2}q^{\fm}_{[c_{-}]}\in Q(\partial^{-}_{\infty}M,[c_{-}])$ respectively.\\\\ We then consider the case of {quasi-Fuchsian half-pipe} manifolds (see Definition \red{\ref{quasifuchsianhalfpipe}}, also\cite{Danciger2013,barbot,Fillastre2019}). These are intermediary geometric structures that arise naturally when we consider smooth transitions between hyperbolic and {anti-de Sitter} structures on $M$ via the Fuchsian locus; the bending laminations on the convex core boundary of the latter being a well studied topic as well. We define an analogous notion for {half-pipe Schwarzians} in this situation (see Definition \red{\ref{definitionofhpsch}}) and prove:
	{\theorem\label{thm1.3} Given any pair of filling measured foliations $\fp,\fm$, there exists a unique quasi-Fuchsian half-pipe manifold such that the horizontal measured foliations of its positive and negative half-pipe Schwarzians are given by $\fp$ and $\fm$ respectively.}

	\subsection{Analogy between bending laminations and measured foliations at infinity}\label{schwatinf}\hfill \break \\
	\hspace*{.3 cm} There are a few points of analogies between the data on the boundary at infinity and that on the boundary of the convex core of a quasi-Fuchsian manifold which makes Theorem \ref{thm1.1} really interesting. We denote the {convex core} of $M$, as $\mathcal{CC}(M)$, as the smallest non-empty convex compact subset contained in $M$ and it is homeomorphic to $S\times [-1,1]$. Call $\partial^{+}\mathcal{CC}(M)$ and $\partial^{-}\mathcal{CC}(M)$ (see \S \ref{introtofuchsian}) as the two boundary components and let the induced metric be called $m_{+}$ and $m_{-}$ respectively. There is a conjecture of Thurston regarding parametrization of quasi-Fuchsian metrics on $M$ uniquely by the data $(m_{+},m_{-})$ (see\cite{Canary2011,sullivan,Labourie1992}). The components $\partial^{\pm}\mathcal{CC}(M)$ moreover carry two {measured geodesic laminations} $\lambda_{+}$ and $\lambda_{-}\in\ML$ where, $\mathcal{ML}(\Sig)$ is the space of measured geodesic laminations on $S$ up to  equivalence (see\cite{bonahon05}). These are called the bending laminations and $\partial^{\pm}\mathcal{CC}(M)$ are bent along leaves of $\lambda_{\pm}$ respectively with the bending angle being given by the transverse measures associated to $\lambda_{\pm}$.\\\hspace*{.5cm} The similarity between the variational formulae for the dual volume $V^{*}_{C}(M)$ of $\mathcal{CC}(M)$ (see \cite{Krasnov2009}) and the renormalised volume $V_{R}$ of $M$ (see \cite{Krasnov2008}) makes Theorem \ref{thm1.1} really interesting as well. Suppose for $0\leq t<\epsilon$, we have a differentiable path of quasi-Fuchsian metrics on $M$ given by $t\mapsto M_{t}$, then the formula for the first-order variation of the renormalised volume is given by (\cite{schlenker20}): 
	\begin{align}\label{renormazlised volume}
		\ddt V_{R}(M_{t}) = -\frac{1}{2}d(\ext(\fp))(\ddt [c^{t}_{+}])	
	\end{align} where $[c^{t}_{+}]$ denotes the variation of the complex structure (up to equivalence) on $\partial^{+}_{\infty}M_{t}$ and for a measured foliation $\f\in \MF$ we have the function $\ext(\f): \T \rightarrow \mathbb{R}$ sending a conformal class $[c]\in \T$ to the extremal length $\ext_{[c]}(\f)$ of the foliation in that conformal class (see \S \S \ref{extremal}). On the other hand, the first order variation of the dual volume, via an application of the Bonahon-Schl{\"a}fli formula is expressed as (\cite{Krasnov2009,Mazzoli2021}):
	\begin{align}\label{dualvolume}
		\ddt V^{*}_{C}(M_{t})=-\frac{1}{2}d(l(\lambda_{+}))(\ddt m^{t}_{+})
	\end{align} where for a measured geodesic lamination $\lambda\in \ML$ we have the function $l(\lambda): \T \rightarrow \mathbb{R}$ sending a hyperbolic metric $m\in \T$ to the length of $\lambda$, denoted as $l_{m}(\lambda)$, measured with respect to this metric and $m^{t}_{+}$ denotes the variation of the induced metric on the convex core boundary under the variation of the quasi-Fuchsian structure. Here, we note that for a given measured foliation $\f$ and measured lamination $\lambda$ the derivatives $d(\ext(\fp)),d(l(\lambda_{+})):T\T \rightarrow \mathbb{R}$ are considered as elements in the cotangent space $T^{*}\T$. Moreover, we also have the upper bound from\red{\cite{Bridgeman2019}} that  $l_{m_{\pm}}(\lambda_{\pm})\leq 6\pi|\chi(\Sig)|$ whereas, from\red{\cite{schlenker20}} we have similar upper bounds on the extremal length $\ext_{[c_{\pm}]}(\f_{\pm})\leq3\pi|\chi(\Sig)|$, where $\chi(S)$ is the Euler characterisitic of $S$. \\ \hspace*{.5 cm}Further, there is a well-studied conjecture of Thurston which asks if the map $
	\mathfrak{B}:\QF \rightarrow \mathcal{ML}(\Sig) \times \mathcal{ML}(\Sig)
	$ sending a quasi-Fuchsian metric $g\in \QF$ to the data $\mathfrak{B}(g):=(\lambda_{+},\lambda_{-})$ of measured bending laminations on the boundary of its convex core, is a homeomorphism onto its image? That is to say, whether quasi-Fuchsian metrics on $M$ can be parametrized by the data of measured bending laminations $(\lambda_{+},\lambda_{-})$ on the boundary of its convex core.
	Although the problem remains open in full generality (see also \cite{Bonahon2004,Lecuire2005,Series,geuritaud} and \cite{adsbending} for the anti-de Sitter case) it can be seen from rather elementary arguments that the image of the map $\beta(\QF)$ is contained in $\mathcal{FML}(\Sig)$, the space of pairs of {filling} measured geodesic laminations on $\Sig$, i.e,  $\lambda_{+}$ and $\lambda_{-}$ always fill $S$ for any quasi-Fuchsian manifold. Using this property Bonahon proves the following theorem  to which we claim our Theorem \red{\ref{thm1.1}} is an analogue of when restricted to the case of measured foliations which are arational: 
	{\theorem \red{\cite{bonahon05}}\label{bonahon}
		There exists an open neighbourhood $V$ of $\mathcal{F}(S)$ in $\mathcal{QF}(S)$, such that $\mathfrak{B}:\QF\rightarrow \ML\times\ML$ is a homeomorphism between $V\setminus \mathcal{F}(S)$ and its image.  Moreover, $V$ can  be chosen so that, $\mathfrak{B}(V \setminus \F)=U$ is an open subset of $\mathcal{FML}(S)$ which intersects
		each ray $(0,\infty)(\lambda_{+},\lambda_{-})$ in an interval $(0,\epsilon_{\lambda_{\pm}})(\lambda_{+},\lambda_{-})$.}\\\\
	\hspace*{.5 cm} A consequence of the theorem above is that the image $\mathfrak{B}(U\setminus\F)$ are pairs of filling measured geodesic laminations $(t\lambda_{+},t\lambda_{-})$, for $t>0$ small enough and clearly, this inspires Theorem \ref{thm1.1}. measured foliations at infinity of $M$ can be thus thought of as a new invariant that provide coordinates for $\QF$ near the Fuchsian locus in a fashion similar to that of measured bending lamination on the boundary of the convex core $\mathcal{CC}(M)$ and we summarise the preceding discussion as Table \ref{tb:1}. We conjecture that our current result can be extended to any pair $(t\fp,t\fm)\in \FMF$ for $t$ small enough.
	\begin{table}[H]\label{table}
		
		\begin{tabular}{|c|c|}
			\hline
			
			On the convex core & On the boundary at infinity \\
			\hline\hline
			Thurston's conjecture on $(m_{+},m_{-})$ & Bers' Simultaneous Uniformisation Theorem \\ \hline 
			Hyperbolic length $l_{m_{\pm}}(\lambda_{\pm})$ & Extremal length  $\ext_{[c_{\pm}]}(\f_{\pm})$
			\\\hline
			$l_{m_{\pm}}(\lambda_{\pm})\leq 6\pi|\chi(S)|$ &   $\ext_{[c_{\pm}]}(\f_{\pm})\leq 3\pi|\chi(S)|$\\ \hline 
			
			Variational formula (\ref{dualvolume}) for $V^{*}_{C}$ & Variational formula (\ref{renormazlised volume}) for $V_{R}$\\  \hline
			
			Theorem \ref{bonahon} & Theorem \red{\ref{thm1.1}} \\
			\hline
		\end{tabular}
		\caption{}
		\label{tb:1}
	\end{table} 
	
	\subsection{Outline}\label{2}
	We prove Theorem \ref{thm1.1} by showing the existence of unique paths in $\QF$ starting from the Fuchsian locus whose measured foliations at infinity are given by $(t\fp,t\fm)\in \mathcal{FMF}_{0}(S)$ for $t>0$ is small enough. Following\cite{bonahon05} this is done essentially by applying an inverse function theorem to the function $\mathfrak{F}: \QF \rightarrow \MF \times \MF$ at the Fuchsian locus and to remove the non-degeneracy of $\mathfrak{F}$ at $\F$, we pass to the blow-up space $\widetilde{\QF}$. To methodize, in \red{\S \ref{minsec}} we establish a necessary condition that infinitesimal deformations of quasi-Fuchsian metrics starting from the Fuchsian locus should satisfy if they have any pair of filling measured foliations $(t\fp,t\fm)$ appearing as their foliation at infinity at first order at $\F$ (Proposition \ref{neccesary}). In \red{\S \ref{pathexists}} we then use this condition to construct small paths $g_{t}$ of quasi-Fuchsian metrics starting from the Fuchsian locus which satisfies $\mathfrak{F}(g_{t})=(t\fp,t\fm)\in \mathcal{FMF}_{0}(S)$ for $0<t<\epsilon_{\fpm}$ where $\epsilon_{\fpm}$ depends on $(\fp,\fm)$ (we don't know how the $\epsilon_{\fpm}$ depends on the pair though). For this we study the sections $q^{\f},q^{-\g}:\T \rightarrow T^{*}\T$ for an arational filling pair $(\f,\g)$. An important step in the proof is to identify the intersection of $[q^{\f}]$ and $[q^{-\g}]$ in the quotient unit bundle $UT^{*}\T$ with a Teichm{\"u}ller geodesic given by the critical point of the function $\ext(t \f)+\ext(\g):\T \rightarrow \mathbb{R}$; this is done in  \red{\S \ref{horifoli}}. In \red{\S\ref{halfpipe}} we interpret our results in half-pipe geometry using tools from the previous sections. We define the notion of half-pipe Schwarzians (see \red{\S \S \ref{hpschwarzandfoli}}) once more by using the paths $\beta_{([c],q)}(t)$ and prove Theorem \red{\ref{thm1.3}}. In \S \ref{halfpipe} we use the results in \S \ref{minsec} to prove our result on half-pipe quasi-Fuchsian manifolds. \S \ref{prelim} contains the necessary preliminaries. \\
	
	\textbf{Acknowledgements}: I heartily thank my advisors Greg McShane, Jean-Marc Schlenker and Andrea Seppi for their encouragements, patience and time. I also thank François Filllastre and anonymous referees for their inputs, and pointing out the errors in the previous version. Finally, this paper wouldn't have been possible without the help from many mathematicians; in particular, I thank Thierry Barbot, Chinmoy Bhattacharjee, Francesco Bonsante, Francis Bonahon, Tommaso Cremaschi, David Dumas, Christian El Emam, François Geuritaud, Gianluca Faraco, Cyril Lecuire, Brice Loustau, Erwan Lanneau, Filippo Mazzoli and Nathaniel Sagman for fruitful discussions.

	\section{Preliminaries}\label{prelim}

	\subsection{Hyperbolic surfaces} \hfill \break 
	
	We will consider $S$ to be a closed surface of genus $g\geq 2$. We call $S$ to be a hyperbolic surface if we have an atlas $(U_{i},\phi_{i})$ on $S$ where $\phi_{i}: U_{i}\rightarrow \mathbb{H}^{2}$ are charts such that at each intersection $U_{i}\cap U_{j}$, the composition $\phi_{i} \circ \phi^{-}_{j}$ are locally restrictions of elements of $PSL_{2}(\mathbb{R})$. An alternative defition can be to say that $S$ is a closed hyperbolic surface if it carries a complete Riemannian metric of constant sectional curvature $-1$. It follows from Gauss-Bonnet theorem that $S$ can carry such a metric only when $g>1$. In such a case, one can also state that $S$ is isometric to the quotient of $\mathbb{H}^{2}$ by $\Gamma$ where $\Gamma$ is a discrete subgroup of $PSL_{2}(\mathbb{R})$.\\
	On the other hand, a complex structure $c$ on $\Sig$ consists of an atlas $\left\lbrace U_{\alpha},\phi_{\alpha}\right\rbrace $ on $\Sig$ where $\phi_{\alpha}:U_{\alpha}\rightarrow \mathbb{C}$ are holomorphic maps and the transition functions $\phi_{i}\circ \phi^{-1}_{j}$ are biholomorphic maps on $\phi_{i}(U_{i}\cap U_{j})$. Given a complex structure $c$, we consider its equivalence class under diffeomorphisms of $S$ isotopic to the identity and denote it as $[c]$. \\
	Consider now a Riemannian metric $g$ on $S$, we can define:
	{\definition A conformal class on a surface $\Sig$ is an equivalence class of Riemannian metrics $[[g]]$, where \begin{align*}
			[[g]]=\left\lbrace e^{2u}g| u: S \rightarrow \mathbb{R}\right\rbrace 
		\end{align*} and $u$ is smooth.} \\ \\
	When $S$ is oriented there is a one-to-one correspondence between equivalence classes of complex structures on $S$ under diffeomorphisms isotopic to the identity and conformal classes on $S$  again up to diffeomorphisms isotopic to the identity. Owing to this we can also view Teichm{\"u}ller space as the space of conformal classes on $S$ (see \cite{bookoftromba}). We will denote both a conformal and complex structure on $S$ as $c$. A Riemannian metric on $S$ in a conformal class has the local expression $g(z)=\rho(z)dzd\bar{z}$ where $\rho(z)\geq 0$ is a smooth function on $\Sig\rightarrow \mathbb{R}_{> 0}$. 
	\\\hspace*{.5 cm}We will recall now an important lemma concerning the change of {Gaussian} or intrinsic curvature $K_{g}$, associated to a Riemannian metric $g$ under change of conformal factor in the same conformal class. See\red{\cite{Krasnov2007}} among others for a reference:
	{\lemma\label{formula} Let $g$ and $g'$ be two Riemannian metrics on $\Sig$ in the same conformal class and let  $u: \Sig \rightarrow \mathbb{R}$ be a function such that $g'=e^{2u}g$. Let $K_{g}$ and $K_{g'}$ be the Gauss curvatures associated to $g$ and $g'$ respectively. Then $K_{g'}=e^{-2u}(-\Delta_{g} u+ K_{g})$, where $\Delta_{g}u$ is the Laplace-Beltrami operator for the metric $g$.}\\\\
	Here we use the convention that $\Delta_{g}$ is negative of the usual analysts Laplacian. That is given the Levi-Civita connection $\nabla$ of $g$ we define $\Delta_g u = -\trace( \hess (u))=-\trace(\nabla \nabla u) $. If we consider $g$ to be a conformal metric, then the hyperbolic metric $m$ in the conformal class of $g$ is given by $m=e^{2u}g$ where $u$ solves:
	{\begin{align*}
			-1=e^{-2u}(-\Delta_{g} u+K_{g})
	\end{align*}}
	So in dimension $2$, corresponding to every conformal class on $S$, one has a unique hyperbolic metric and also an equivalence class of complex structures. \\ 
	
	\subsection{Fuchsian and quasi-Fuchsian $3$-manifolds}\label{introtofuchsian} \hfill \break \\
	We will consider $3$ manifolds $M$ here which are quotient $\mathbb{H}^{3}/\Gamma$ where $\Gamma$ is a discrete subgroup of $PSL_{2}(\mathbb{C})$. Associated to the action of $\Gamma$ on $\mathbb{H}^{3}$ one has the limit set $\Lambda_{\Gamma}$ which is the set of accumulation points of orbit of $\Gamma$, and it can so shown that it is a subset of $\partial_{\infty}\mathbb{H}^{3}$. When $\Gamma$ is a discrete subgroup of $PSL_{2}(\mathbb{R})\subset PSL_{2}(\mathbb{C})$, we call $M$ to be a Fuchsian manifold. In this case $\Lambda_{\Gamma}$ is a circle on $\partial_{\infty}\mathbb{H}^{3}$. This can be seen as the boundary of the totally geodesic copy of $\mathbb{H}^{2}\subset \mathbb{H}^{3}$ preserved by the action. \\
	We call a $3$ manifold $M$ to be quasi-Fuchsian if the $\Gamma < PSL_{2}(\mathbb{C})$ is such that $\Lambda_{\Gamma}$ is a quasi-circle. To define concretely, one has a quasi-conformal map $\phi: \mathbb{C}P^{1}\rightarrow \mathbb{C}P^{1}$ and a Fuchsian subgroup $\Gamma_{0}$ such that $\Gamma := \phi^{-1}\circ \Gamma_{0} \circ \phi$. To simplify, we can consider $\Gamma$ to be quasi-Fuchsian if $\Lambda_{\Gamma}$ is a Jordan curve on $\partial_{\infty}\mathbb{H}^{3}$. \\ 
	One more charateristic difference between Fuchsian and quasi-Fuchsian manifolds is by the different geometry of their convex hulls. Given the limit set $\Lambda_{\Gamma}$ as the setting above, we can consider its convex hull in $\mathbb{H}^{3}$. In the Fuchsian case we recover the totally geodesic copy of $\mathbb{H}^{2}$ preserved by the action as $\Lambda_{\Gamma}$ is a circle. In the quasi-Fuchsian case we have that the convex hull is a closed, convex region which upon taking quotient becomes the convex core $\mathcal{CC}(M)$. The convex core is thus the smallest-non empty convex submanifold contained in the quasi-Fuchsian manifold and it comes with two boundary components since it is $S\times [-1,1]$. 
	\subsection{Teichm{\"u}ller space} \label{teichspace}\hfill \break \\
	We will now briefly introduce the Teichm{\"u}ller space $\T$ associated to a closed surface of genus $g\geq 2$.
	{\definition The Teichm{\"u}ller space $\T$ is the space equivalence classes of complex structures on $S$ under diffeomorphisms isotopic to the identity.}\\\\
	Owing to the presence of an unique hyperbolic metric in every conformal class or class of complex structure we can alternately define $\T$ as the space of equivalence classes of hyperbolic metrics on $S$ up to diffeomorphisms isotopic to the identity.\\
	Now recall that the boundary at infinity $\partial^{\pm}_{\infty}M$ are two copies of $S$ which carry their respective complex structures. It so happens that when $M$ is Fuchian, the complex structures/hyperbolic metrics are identical and thus if one defines $\F$ as the equivalence class of Fuchsian metrics on $M$ up to diffeomorphism isotopic to the identity, then $\F\cong \T$. \\
	Now, given two conformal structures $c,c'$ we can define the notion of quasi-conformal map between them.
	{\definition Let $\Omega,\Omega'\subset\mathbb{C}$ be two domains and $\phi: \Omega\rightarrow \Omega'$ be a homeomorphism with continuous partial derivatives with respect to $z$ and $\bar{z}$. We denote $\phi$ to be $k$-quasiconformal if 
		\begin{align*}
			k_{\phi}(z)=\frac{|\frac{\partial \phi}{\partial \bar{z}}(z)|+|\frac{\partial \phi}{\partial z}(z)|}{|\frac{\partial \phi}{\partial \bar{z}}(z)|-|\frac{\partial \phi}{\partial z}(z)|}\leq k
		\end{align*}
		for almost every $z\in \Omega$.}\\\\ 
	We will denote the $\frac{\partial \phi}{\partial w}$ as $f_{w}$. The {Beltrami differential} $\mu$ associated to the map $\phi$ is defined as the ratio 
	\begin{align*}
		\mu_{\phi} = \frac{\phi_{\bar{z}}}{\phi_{z}}
	\end{align*}
	which is defined almost everywhere, is measurable and satisfies $||\mu||_{\infty}<1$. Equivalently $k_{\phi}$ has the expression $\frac{1+|\mu|}{1-|\mu|}$, which is bounded above by $k=\frac{1+||\mu||_{\infty}}{1-||\mu||_{\infty}} $ and $k_{\phi}$ is called the eccentricity coefficient of $\phi$. We also note that a Beltrami differential on $(S,[c])$ is a tensor of type $(-1,1)$. We will use this to define the notion of holomorphic quadratic differentials on $S$ and how it relates to the tangent and cotangent space of $\T$.
	\subsection{Holomorphic quadratic differentials}\label{HQDintro} \hfill \break \\
	A {holomorphic quadratic differential} $q$ on $(S, c)$ is a tensor of type $(2,0)$ which in local coordinates can be written as $f(z)dz^{2}$, where $f$ is a holomorphic function. The space of holomorphic quadratic differential denoted as $Q(S)$ forms a bundle over $\mathcal{T}(S)$, where $\T$ is seen as the space of complex structures on $\Sig$ up to  diffeomorphisms isotopic to the identity. The fiber over an equivalence class $[c] \in \T$, which is denoted as $Q(S, [c])$ is a vector space that can be shown to have real dimension $6g-6$ (by, for example, the {Riemann-Roch formula}). Moreover, a holomorphic quadratic differential $q$ has zeroes on $\Sig$ the degree of which is defined in terms of the degree of the zero of the local Taylor expansion of $q$. To be precise, if $q$ has a zero of order $k$ at a point $p\in S$ then this means that for all chart centred at $p$ on $S$, $q$ has the local expression $f(z)z^{k}dz^{2}$ for some holomorphic function $f(z)$ such that $f(0)\neq 0$.. Moreover, it follows from, for example the Riemann-Roch theorem, that the sum of the degrees of all zeroes of $q$ on $\Sig$ is $4g-4$.  The space $Q(S)$ further carries a natural stratification depending on the order of the zeroes of $q$. Please consult, for example \cite{Hubbard} for references on this topic. 
	{\definition Let $\overline{k}$ be a $n$-tuple of integers $(k_{1},k_{2},\dots,k_{n})$ such that $\sum_{i}^{n}k_{i}=4g-4$. The stratum $Q^{\overline{k}}(S)$ is the set of holomorphic quadratic differentials $q$ such that the degrees of the zeroes of $q$ are given by the $k_{i}$. We say that $q$ is generic, if $k_{i}=1$ for all $i$.} \\\\
	Holomorphic quadratic differentials with only $4g-4$ simple zeroes are termed as generic quadratic differentials and it is known from\cite{Douady1975} that they form a dense open subset of $Q(S)$ which will be denoted as $Q_{0}(S)$.\\\\
	Notice that the product of a Beltrami differential and a holomorphic quadratic differential gives us a $(1,1)$ tensor. Let $B(S,c)$ denote the vector space of measurable Beltrami differentials on $(S,c)$ where an element is expressed locally as $\mu = b(z)\frac{d\bar z}{d z}$ . From here we have a natural complex pairing between $\mu \in B(S,c)$ and $q \in Q(S,c)$ as: 
	\begin{align*}
		\left\langle q,\mu \right\rangle = \int_{(S,c)} q\mu dz d\bar{z} 
	\end{align*}
	It follows as a consequence, see \cite{Hubbard}, that: 
	{\proposition \label{keylemma}There is an isomorphism of vector spaces between 
		\begin{align*}
			T_{[c]}\T \cong B(S,c)/Q(S,c)^{\perp} {\quad and \quad} T^{*}_{[c]}\T \cong Q(S,c)
		\end{align*}
		where $Q(S,c)^{\perp}= \left\lbrace \mu \in B(S,c)|\left\langle \mu,q\right\rangle  =0, \forall q \in Q(S,c)\right\rbrace$.}\\\\
	This allows us to define the Weil-Petersson metric on $T^{*}_{[c]}\T$ as: 
	\begin{align*}
		\left\langle q_{1},q_{2}\right\rangle_{WP} = \int_{S} \frac{f_{1}(z)\overline{f_{2}(z)}}{\rho(z)}dz \overline{dz}
	\end{align*}
	where the hyperbolic metric in the class $[c]$ has the expression $\rho(z)|dz|^{2}$ and $q_{i}=f_{i}(z)dz^{2}$ for $i=1,2$ are two holomorphic quadratic differentials in $Q(S,c)$. This also induces an inner product on the tangent space $T_{[c]}\T$ by duality. The Weil-Petersson metric gives $\T$ with the structure of negatively curved Riemannian manifold. On the other hand, the $L^{1}$-norm on the cotangent space $T^{*}_{[c]}\T$ is defined as:
	\begin{align*}
		||q||_{1} = \int_{(S,c)}|f(z)| dz \wedge \overline{dz}.
	\end{align*}
	This induces a Teichm{\"u}ller norm on $T\T$ via the duality between Beltrami differentials and holomorphic quadratic differentials. One way to express the associated metric, called the Teichm{\"u}ller metric $d_{\T}$, is:
	\begin{align}\label{teichmullergeo}
		d_{\T}([c],[c']):=\frac{1}{2}\inf\left\lbrace\log k_{\phi}| \text{ }\phi:(\Sig,[c])\rightarrow (\Sig,[c'])\text{ quasiconformal isotopic to the identity} \right\rbrace 
	\end{align} One can consult, for example \cite{carlos}, for further details in this topic.

	\subsection{Measured foliations on $S$ and the space $\MF$}\label{defmf}\hfill\break \\
	Following\red{\cite{hubbardmasur}} we define:
	{\definition\label{deffoliations} A smooth measured foliation $\f$ on $\Sig$ with singularities $\left\lbrace p_{1},\dots,p_{n}\right\rbrace $ of
		order $\left\lbrace k_{1},\dots,k_{n}\right\rbrace $ (respectively) is given by an open covering ${U_{i}}$ of $\Sig\setminus\left\lbrace p_{1},\dots,p_{n}\right\rbrace $ and
		open sets $\left\lbrace V_{1},\dots,V_{n}\right\rbrace $ around $\left\lbrace p_{1},\dots,p_{n}\right\rbrace $ (respectively) along with smooth non-vanishing real valued $1$-forms $d\phi_{i}$
		defined on $U_{i}$ such that:
		\begin{itemize}
			\item $d\phi_{i}=\pm d\phi_{j}$ on $U_{i}\cap U_{j}$
			\item around each $p_{l}$ there is an open neighbourhood $V_{l}$ and a chart $(x_{1},x_{2}):V_{l}\rightarrow \mathbb{R}^{2}$ such that $d\phi_{i} = \mathfrak{I}(z^{\frac{k_{l}}{2}+1} dz)$ on $U_{i}\cap V_{j}$ where $z=x_{1}+\textbf{i}x_{2}$.
	\end{itemize}}

	Immersed lines on $\Sig$ along which $d\phi_{i}$ vanish give a
	foliation $\f$ on $\Sig \setminus (p_{1},\dots,p_{n})$ and we have a $(k_{j}+2)$ pronged singularity at $p_{j}$. Given an arc $\gamma$ on $\Sig$ which avoids the zeroes $(p_{1},p_{2},\dots,p_{n})$, a measured foliation $\f$ associates a transverse measure to $\gamma$ defined as $\mu_{\f}(\gamma)= \int_{\gamma}|d\phi|$, where $|d\phi|$ restricted on each $U_{i}$ is given by $|d\phi_{i}|$. \\\hspace*{.5 cm}
	This measure is invariant under isotopies that maintain the same end points of $\gamma$ and the transversality of the intersection of $\gamma$ with the given foliation. That is, if $\gamma'$ is isotopic to $\gamma$ with the same end points and maintaining the transversality at every time, then $\mu_{\f}(\gamma)=\mu_{\f}(\gamma')$. So, given a class $[\gamma]$, we define:
	{\definition {The intersection number} of $\f$ with a isotopy class of closed curves $[\gamma]$ avoiding the singularities of $\f$ is defined $i([\gamma],\f)=\inf_{\gamma \in [\gamma]}i(\gamma,\f)=\inf_{\gamma \in [\gamma]}\mu_{\f}(\gamma)$, where the infimum is taken over all $\gamma \in [\gamma]$. }\\\\
	As we can see that the intersection number defines a function from the set of closed curves up to isotopy on $\Sig$ to $\mathbb{R}_{>0}$ and we define following, for example\red{\cite{AST_1979__66-67_}}:
	{\definition 
		Two measured foliations $\f$ and $\g$ on $\Sig$ are said to be equivalent if they define the same intersection number. The space of equivalence classes of measured foliations on $\Sig$ will be denoted as $\MF$. 
	}\\
	
	The space $\MF$ can also be defined through a topological equivalence of two foliations which comes via {Whitehead moves} but we do not elaborate on that. Also it follows from Proposition {$2.1$} in\red{\cite{hubbardmasur}} that a measured foliation $\f$ can have $4g-4$ pronged-singularities counted up to multiplicity. \\\hspace*{.5 cm}
	There is also an action of $\mathbb{R}_{>0}$ on $\MF$ defined as $t.\f\mapsto t\f$ where the latter denotes the measured foliation obtained by multiplying $t$ by the $1$-forms $d\phi_{i}$ which give us the measured foliation $\f$ as par Definition \ref{deffoliations}. 
	The space $(\MF\setminus 0)/\mathbb{R}_{>0}$ is called the space of projectivised measured foliations, denoted as $P\MF$ and it is identified with the {Thurston boundary} of $\T$. We will denote by $[\f]$ as the equivalence class of $\f$ in $P\MF$. Interested readers can consult, for example\red{\cite{AST_1979__66-67_}}, for further details in this topic.  \\
	We also introduce the notion of arational measured foliations here as: 
	{\definition \label{arationaldef} A measured foliation is said to be arational if all its singularities have $3$ prongs and if there are no leaves of the foliation joining the singularities.}\\ 
	
	We call leaves of the foliations joining prongs of the singularities as saddle connections.

	\subsection{Minimal surfaces in hyperbolic $3$ manifolds}\label{minimalsurface}\hfill\break\\
	Let $i: \Sig \rightarrow M$ be an immersion of $\Sig$ into $M$. Associated to this, we have the data of the {first fundamental form } denoted as $I$ on $\Sig$ which is the induced metric on $\Sig$ inherited from the ambient space and the {second fundamental form} denoted a $\II$ which is a symmetric bilinear form on the tangent bundle $T\Sig$. Associated to the couple $(I,\II)$ we have a unique self-adjoint operator $B$, called the {shape operator}, which satisfies the relation \begin{align*}
		\II(x,y)=I(Bx,y)=I(x,By)
	\end{align*} where $x,y \in T_{p}S, \forall p \in S$. We can also define another quantity associated to the immersion called the third fundamental form $\III$ defined as $\III(x,y)=I(Bx,By)=I(B^{2}x,y)$ for $x,y,B$ as before.\\\hspace*{.5 cm}
	The eigenvalues of $B$ gives us the principal curvatures associated to the immersion. Since minimal immersions are those immersions for which the mean curvature of $\Sig$ is zero we can define as well that: 
	{\definition An immersion is minimal if and only if $B$ is traceless. }
	\\
	
	On the other hand, assume we are given a smooth Riemannian metric $g$ on $\Sig$ and a symmetric bilinear form $h$ on $T\Sig$. The couple $(g,h)$ will be associated to the data of an immersion  of $\Sig$ if it satisfies the following:
	
	\begin{enumerate}
		
		\item The Codazzi equation, $d^{\nabla}h=0$, where $\nabla$ is the Levi-Civita connection of $g$.
		\item The Gauss equation, $K_{I}=-1+\deter_{g}(h)$, where $K_{g}$ denotes the Gaussian or intrinsic curvature of the metric $g$.
	\end{enumerate}
	Now we define the notion of an {almost}-Fuchsian manifold as: 
	
	{\definition A quasi-Fuchsian hyperbolic $3$-manifold $M \cong \Sig \times \mathbb{R}$ is called {almost-Fuchsian} if it contains a closed minimal surface homeomorphic to $\Sig$ with principal curvatures in $(-1,1)$.}\\\\
	Given an immersed surface in $M$, an outcome of the Codazzi equation is that, the traceless part of the second fundamental form $(\II)_{0}$ is equal to the real part of a holomorphic quadratic differential $q$. 
	Now, following Corollary $2.9$ of\red{\cite{Krasnov2007}} we can show that if $M$ is an almost-fuchsian manifold then the closed minimal surface it contains is unique. This allows us to parametrise almost-Fuchsian hyperbolic manifolds by an open subset $\Omega$ in $T^{*}\T$, following \cite{Krasnov2007}, Theorem $2.12$:
	{\theorem\label{methodandoutline} There exists an open subset $\Omega \subset T^{*}\T$ such that we have the following bijection:
		\begin{enumerate}
			\item Given a point $([c],q)\in \Omega$, there exists a unique almost-Fuchsian metric on $M$ such that the unique minimal surface has first fundamental form conformal to $[c]$ and the second fundamental form $\II$ is $\mathfrak{R}(q)$.
			\item Given a almost-Fuchsian metric $g \in \AF$ on $M$, the induced metric and  second fundamental form of its unique minimal surface are specified by a point in $\Omega$. 
	\end{enumerate}} 
	
	\subsection{Schwarzians at infinity of quasi-Fuchsian manifolds}\label{introductionary}\hfill \break \\
	Recall that the boundary at infinity $\partial_{\infty}\mathbb{H}^{3}$ is identified with the complex projective space $\mathbb{C}P^{1}$. The components $\partial^{+}_{\infty}M$ and $\partial^{-}_{\infty}M$ are the quotient of domains in $\mathbb{C}P^{1}$ under the action of $\Gamma < PSL_{2}(\mathbb{C})$ and so they carry canonical $\mathbb{C}P^{1}$-structure, that is we have an open covering of $\Sig$ by an atlas $(U_{\alpha},\phi_{\alpha})$ such that $\phi_{\alpha}:U_{\alpha}\rightarrow \mathbb{C}P^{1}$ are charts to open domains in $\mathbb{C}P^{1}$ and in the overlap  $U_{i}\cap U_{j}$  of two charts the change of coordinate map $\phi_{i}\circ \phi^{-1}_{j}$ is locally a restriction of a M{\"o}bius transformations. Denote the space of equivalence classes of $\mathbb{C}P^{1}$-structures on $\Sig$ under diffeomorphisms isotopic to the identity as $\mathcal{CP}(\Sig)$.  \\\hspace*{.5 cm}Now, given a $\mathbb{C}P^{1}$-structure on $\Sig$ we have an underlying complex structure as M{\"o}bius transformations are biholomorphisms and for $\partial^{\pm}_{\infty}M$ it is precisely $[c_{\pm}]$ up to equivalence. This gives us a natural forgetful map $\mathcal{CP}(\Sig)\rightarrow \T$ mapping a $\mathbb{C}P^{1}$-structures to the underlying complex one. Now by the  Uniformisation theorem any complex structure $c$ on $\Sig$ arises as the quotient of the action of some discrete subgroup  $\Gamma_{c}$ of $PSL_{2}(\mathbb{R})$ on $\mathbb{H}^{2}$. As $\Gamma_{c}<PSL_{2}(\mathbb{C})$ as well and $\mathbb{H}^{2}$ can be seen as the unit disc $\Delta\subset \mathbb{C}P^{1}$, we have a canonical $\mathbb{C}P^{1}$-structure associated to a complex structure which we call the {standard Fuchsian complex projective structure} and this gives us a continuous section $\T \rightarrow \mathcal{CP}(\Sig)$. 
	The Schwarzian derivative yields a parametrisation of the fibers  of the forgetful map $\mathcal{CP}(\Sig)\rightarrow \T$. In general, given a domain $\Omega\subset \mathbb{C}$, 
	the Schwarzian derivative of a locally injective holomorphic map $u:\Omega \rightarrow \mathbb{C}$ is a holomorphic quadratic differential defined as:
	\begin{align*}
		\sigma(u)= ((\frac{u''}{u'})^{'}-\frac{1}{2}(\frac{u''}{u'})^{2})dz^{2}
	\end{align*} One way to obtain the expression on the right hand side above is to consider the unique M{\"o}bius transformation
	$M_{u}$ which matches with $u$ up to second order derivative. The expression above is precisely the difference of
	the third order terms in the local Taylor series expansion of $u$ and $M_{u}$ (see Proposition $6.3.3$ of\cite{Hubbard}). Further, they have two remarkable properties:
	\begin{itemize}
		\item For two locally injective holomorphic maps $u,v: \Omega \rightarrow\mathbb{C}$ we have $\sigma(u\circ v)= v^{*}\sigma(u)+\sigma(v)$
		\item $\sigma(A)=0$ if and only if $A$ is a M{\"o}bius transformation.
	\end{itemize} In particular, there is a unique map defined up to right conjugation by M{\"o}bius transformations between a given complex projective structure on $S$ and the standard
	Fuchsian one which is holomorphic with respect to the underlying canonical complex structure and by virtue of
	the properties above, the Schwarzian derivative for this holomorphic map can be defined in a chart independent
	way. This is called the Schwarzian parametrisation of a complex projective structure on $S$ with respect to
	the standard Fuchsian complex projective structure (see \cite{Dumas,schlenker20}). When $M \in \QF$, the components of $\partial_{\infty}\mathbb{H}^{3}\setminus \Lambda_{\Gamma}$  have non-trivial Schwarzian derivatives
	associated to them by construction which descend to two holomorphic quadratic differentials on $\partial^{\pm}_{\infty}M$ upon
	taking quotients. So we define: {\definition  The Schwarzians at infinity $\sigma_{+}$ and $\sigma_{-}$ are the holomorphic quadratic differentials obtained on $(\partial^{+}_{\infty}M,[c_{+}])$ and $(\partial^{-}_{\infty}M,[c_{-}])$ by the Schwarzian parametrisation of the $\mathbb{C}P^{1}$ structures on $\partial^{\pm}_{\infty}M$ with respect to the corresponding standard Fuchsian complex projective structure.}\\\\
	Also note that when $M\in\F$ due to the second property above the Schwarzians at infinity are zero. What is more important to us from this discussion is that due to this we get another parametrisation of quasi-Fuchsian manifolds by $Q(S)$ by considering the Schwarzian derivative and complex structure appearing at one end at infinity. So we can therefore construct a well-defined map
	\begin{equation}\label{eq:schwarzian map}
		\mathfrak{A}:\QF\to Q(\partial^{\pm}_{\infty}M)~.
	\end{equation}
	Here $Q(S)$ denotes the bundle of holomorphic quadratic differentials over $\T$, whose fiber over a point $(S,[h])$ coincides with the vector space $H^0((\Sigma,h),K^2)$, where $K$ denotes the canonical divisor of $(S,[h])$. Consequently, the space $Q(S)$ is a complex manifold of dimension $3g-3$, where $g$ denotes the genus of $S$. In fact, the map $\mathfrak{A} $ turns out to be injective, and $\QF$ and $Q(\partial^{+}_{\infty}M)$ being manifolds of the same real dimension, the invariance of domain theorem implies that its image is an open subset of $Q(\partial^{+}_{\infty}M)$ (see \S $2.1.3$ of \cite{cmcdifian}).\\

	\section{Measured foliations realised by holomorphic quadratic differentials.}\label{horifoli}\hfill\break
	Now, given $q\in Q(S,c)$, away from its zeroes we can always perform a local change of coordinates $z\mapsto w:=\int\sqrt{q}$ on $S$ such that $q=f(z)dz^{2}$ has the local expression $dw^{2}$ with respect to this coordinate. If we write $w=w_{1}+\textbf{i}w_{2}$ then the holomorphic quadratic differential $dw^{2}$ canonically equips $\mathbb{C}$ with two measured foliations: 
	\begin{itemize}
		\item The horizontal measured foliation, which are immersed lines given by  $w_{2}=\text{const.}$, i.e the horizontal lines of $\mathbb{C}$. Its transverse measure being given by $|\mathfrak{I}\sqrt{dw^{2}}|=|dw_{2}|$.
		\item The vertical measured foliation, which are immersed lines along which $w_{1}=\text{const.}$, i.e the vertical lines of $\mathbb{C}$.
		Its transverse measure being given by $|\mathfrak{R}\sqrt{dw^{2}}|=|dw_{1}|$.
	\end{itemize} 
	Moreover, notice that the horizontal measured foliations (resp. vertical measured foliations) of quadratic differential $|dw^{2}|$ gives us all the horizontal lines (resp. vertical lines) on $\mathbb{C}$, thus inspiring the nomenclature. So we define: 
	{\definition The horizontal measured foliation  $\mathrm{hor}_{[c]}(q)$ (resp. vertical measured foliation $\mathrm{ver}_{[c]}(q)$) of $q$ on $(S,c)$ is a smooth singular measured foliation, with singularities at the zeroes of $q$, which is obtained locally by pulling back the horizontal measured foliations (resp. vertical measured foliation) of $dw^{2}$ under the change of coordinate $z\mapsto w:= \int \sqrt{q}$ defined above. The transverse measure for the horizontal measured foliation (resp. vertical measured foliation) is given by $|\mathfrak{I}\sqrt{q}|$ (resp. $|\mathfrak{R}\sqrt{q}|$).} \\\\If the measured foliation $\f$ is realised by a holomorphic quadratic differential $q$ then $\f$ has a prong of order $k+2$ at the point where $q$ has a zero of order $k$. Also, if $q$ is expressed as $dw^{2}$ in local coordinates  then $-q$ is nothing but the differential $-dw^{2}$, whose horizontal (resp. vertical) foliations are given by the vertical lines (resp. horizontal lines) on $\mathbb{C}$. We thus have the simple but important remark:
	
	{\remark\label{remarkk} $\mathsf{hor}_{[c]}(q)$ is measure equivalent to $\mathsf{ver}_{[c]}(-q)$ in $\MF$ for all holomorphic quadratic differential $q\in Q(S,c)$. }\\\\
	Now we recall a well-known theorem of Teichm{\"u}ller that enables us to interpret quasi-conformal deformations in terms of measured foliations (see for example\red{\cite{carlos}}):
	{\theorem Given two conformal classes $[c]$ and $[c']$ in $\T$, there exists an unique quasi-conformal map $\phi:[c]\rightarrow [c']$ with minimal eccentricity coefficient among all quasi-conformal maps from $[c]$ to $[c']$. The associated Beltrami differential $\mu$ is of the form $k\frac{ |q|}{q}$ for some unique holomorphic quadratic differential $q\in T^{*}_{[c]}\T$ with $||q||_{1}=1$ and for some $k\in [0,1)$. The quadratic differential $q$ is denoted as the initial quadratic differential of the map. There is a quadratic differential $q'\in T^{*}_{[c']}\T$ denoted as the terminal quadratic differential with the property that the map $\phi$ takes zeroes of $q$ to zeroes of $q'$ of the same order. In the natural local coordinates $z = x + \textbf{i}y$ of $q$ in the complement of its zeroes, and the natural coordinates $w = x' + \textbf{i}y'$ for $q'$, we have:
		\begin{align*}
			x' = \sqrt{k}x, y' = (1/ \sqrt{k})y
		\end{align*}.}
	
	By virtue of the above, the metric $d_{\T}$ defined in Equation (\ref{teichmullergeo}) is a metric on $\T$. Further, we also have that given $(S,c)$, a quadratic differential $q \in T^{*}_{[c]}\T$ with $||q||_{1}=1$ and $t \geq 0$, there is a conformal class $(S,[c_{t}])$, and a unique extremal map $\phi_{t} : (S,c) \rightarrow (S,[c_{t}])$ such that:
	\begin{align*}
		d_{\T}([c],[c_{t}])=\frac{1}{2}\log k_{\phi_{t}}
	\end{align*}
	Choosing $\log k_{\phi_{t}}=2t$ gives us that the image of the map $\mathbb{R}_{> 0} \rightarrow T^{*}\T\cong Q(S)$ with $t\mapsto [c_{t}]$ is a properly embedded geodesic line in $\T$ with respect to the metric $d_{\T}$ which is called the {Teichm{\"u}ller geodesic} with initial quadratic differential $q$. This also makes the metric complete (see, for example \cite{carlos} for more further details). 
	
	\subsection{The sections $q^{\f}$ and $q^{-\f}$}\label{HMsection} \hfill\break\\
	Define now a map $\mathsf{hor}:Q(\Sig,c)\rightarrow \MF$ which sends a holomorphic quadratic differential $q$ to its horizontal measured foliation $\mathsf{hor}_{c}(q)$ and we consider its image in $\MF$. Then we have from\red{\cite{hubbardmasur,wolf}}:
	{\theorem\label{hubbardmasur} The map $\mathsf{hor}_{c}: Q(S, c) \rightarrow \mathcal{MF}(S)$ is a homeomorphism}.
	
	{\remark Given any measured foliation $\f$ on Riemann surface, we may not find a holomorphic quadratic differential realising it as its horizontal measured foliation, for example notice the example in  \S {$2$} of Chapter {$II$} of \red{\cite{hubbardmasur}}. However as noted in the paper, this issue can be taken care of as according to  Proposition {$2.2$} of\red{\cite{hubbardmasur}} as in the equivalence class of $\f$ in $\MF$ there exists a representative that can be realised by a holomorphic quadratic differential.}\\
	
	We will consider the inverse of this map for our purpose which will provide sections of $T^{*}\T$ for a fixed foliation $\f$. This we define as follows: 
	{\definition For a given equivalence class of foliation $\f \in \MF$, define 
		\begin{align*}
			q^{\f}:\T \rightarrow T^{*}\T
		\end{align*} to be the map, which associates to each equivalence class of complex structure $[c]$ on $\Sig$, the unique holomorphic quadratic differential $q^{\f}_{c}$ such that $\mathsf{hor}_{c}(q^{\f}_{c})$ is measure equivalent to $\f$. } \\
	
	We will denote the holomorphic quadratic differential associated to $[c]$ as $q^{\f}_{[c]}$. In fact, the theorem of Hubbard-Masur holds true if we consider vertical measured foliations instead of horizontal ones and thus we can consider the map 
	\begin{align*}
		q^{-\f}: \T \rightarrow T^{*}\T
	\end{align*} which associates to a complex structure $[c]$ on $\Sig$, the unique holomorphic quadratic differential $q^{-\f}_{[c]}$, such that vertical measured foliation of $q^{-\f}_{[c]}$ on $(\Sig,[c])$ is measure equivalent to $\f$. We can thus reformulate Remark \red{\ref{remarkk}} as: 
	
	{\remark\label{lelele} For a given measured foliation $\f \in \MF$, $q^{-\f}_{[c]}=-q^{\f}_{[c]}$ for any $[c]$ in $\T$.}\\\\
	These sections are $C^{0}$, in particular, we don't know if it is $C^{1}$. It also follows from a result of Masur in\red{\cite{Masur1995}} that when $q$ is generic, then the sections $q^{+\f},q^{-\f}: \T \rightarrow T^{*}\T$ are real-analytic. In fact we can state:
	{\lemma\label{arationaliff} A measured foliation is arational if and only if the holomorphic quadratic differential realising it as the horizontal measured foliation over each point in $\T$ is generic. Moreover in this case the map $q^{\f}: \T\rightarrow T^{*}\T$ is smooth.}\\	
	
	One implication is obvious as if the quadratic differential is generic then the measured foliation it realises has three prongs at each zero. The other side can be seen easily as if $F$ is arational then any Whitehead equivalent measure foliation is isotopic to $F$ (since by definition there are no saddle connections to collapse). We also denote the subset of arational measured foliations as $\mathcal{MF}_{0}(S)$ and note that this is a dense subset of $\MF$ as well. 
	\subsection{Filling measured foliations}\label{gardinermasur}\hfill \break \\
	First we recall that:
	{\definition\label{fillupdef}
		A pair of measured foliations $(\mathsf{F},\mathsf{G})$ is said to {fill $S$} if for any measured foliation $\mathsf{H}\in \MF $ on $S$ we have, 
		\begin{align*}
			i(\mathsf{H},\f)+i(\mathsf{H},\g)>0
		\end{align*}.
	} 
	
	Recall that we denote the space of equivalence classes of pairs of filling foliations as $\FMF$. Notice that the pair $(\mathsf{hor}(q),\mathsf{ver}(q))$ are transverse and can be shown to satisfy the topological property of {filling up $S$}, by the following Lemma {$2.10$} of\red{\cite{grenobleschool}}:
	
	{\lemma\label{2.4} Given a holomorphic quadratic differential $q$ on a Riemann surface $(\Sig,[c])$, the pair $(\mathsf{hor}_{[c]}(q),\mathsf{ver}_{[c]}(q))$ fill $\Sig$.}\\
	
	Given a pair $(\f,\g)\in \FMF$ we can thus ask whether under a fixed complex structure up to equivalence, a pair $(\mathsf{F},\mathsf{G})$ can be realized as the horizontal and vertical measured foliation of the same holomorphic quadratic differential $q$. The answers are affirmative and can be summarized as: 
	
	{\theorem(\red{\cite{Gardiner1991,Wentworth2007})\label{GM}
		A pair $(\mathsf{F},\mathsf{G})$ of measured foliations on $S$ is filling if and only
		there is a complex structure $c$ and a holomorphic quadratic differential $q \in (S, c)$ such that $(\mathsf{F},\mathsf{G})$ are respectively measure equivalent to the vertical and horizontal
		foliations of $q$. Moreover, the class $[c]$ up to  diffeomorphism isotopic to the identity is determined uniquely
		and for each $c \in [c]$ the quadratic differential $q$ realising the filling pair $(\f,\g)$ is also unique. }\\ \\
	So, for a pair $(\f,\g)$ that fill we have: 
	{\corollary\label{gola} The sections $q^{\f}$ and $q^{-\g}$ intersect uniquely in $T^{*}\T$ at the point $([c],q)$ determined by Theorem \red{\ref{GM}}. Moreover $\f=\mathsf{hor}_{[c]}(q)$ and $\g=\mathsf{hor}_{[c]}(-q)$ in $\MF$ where $q=q^{\f}=q^{-\g}=-q^{\g}$.}
	\subsection{Extremal lengths of measured foliations}\label{extremal}\hfill\break\\
	Given a simple closed curve $\gamma$ on $(\Sig,c)$ we define its extremal length as \begin{align*}
		\ext_{c}(\gamma)= \sup_{[[g]]}\frac{l^{2}_{g}(\gamma)}{\text{Area}(g)}
	\end{align*}
	where $l_{g}(\gamma)$ is the length computed with respect to $g$ and the supremum is taken over all Riemannian metrics in the conformal class $c$. This definition of extremal length on closed curve extends to that of a measured foliation (see\red{\cite{Kerckhoff1980}}) where extremal length of a foliation $\f$ defines a continuous function on $\T$ 
	\begin{align*}
		\ext(\f):\T\rightarrow \mathbb{R}\\
		[c]\mapsto \ext_{[c]}(\f)
	\end{align*}
	where $\ext_{[c]}(t\f)=t^{2}\ext_{[c]}(\f)$ for $t>0$. Using the sections $q^{\f}_{[c]}$ we can also express this as (see\red{\cite{hubbardmasur}}):
	{\lemma\label{babagola} For $[c]\in \T$, the extremal length of $\f\in\MF$ is given by \begin{align*}
			\ext_{[c]}(\f)=\int_{(S,c)}|q^{\f}_{[c]}|dz\wedge d\bar{z}=||q^{\f}_{[c]}||_{1}=||q^{-\f}_{[c]}||_{1}
	\end{align*} }
	Here $||q||_{1}=\int_{\Sig}|\phi(z)|dz \wedge d\bar z$ where $q=\phi(z)dz^{2}$. Another simple observation that follows from this is:

	{\corollary\label{gorombhat} $\ext_{[c]}(\f)=\ext_{[c]}(\g)$ where $[c]\in \T$ is determined by Theorem \red{\ref{GM}}. }
	{\proof As $q=q^{\f}_{[c]}=q^{-\g}_{[c]}$ at $[c]$ where the sections $q^{\f}$ and $q^{-\g}$ intersect (see Remark \red{\ref{gola}}), we get the result using Lemma \red{\ref{babagola}}. \qed}\\\\
	We also have a well-known variational formula for extremal lengths originally due to Gardiner (see \cite{gardiner1984}) which states:
	{\lemma \label{gardform}Let $[c_{t}]$ for $0\leq t< \epsilon $ be a smooth $1$-parameter family of conformal classes and $\f$ be a smooth measured foliation in $\MF$ then 
		\begin{align*}
			(d \ext_{[c_{0}]}\f) (\mu)=\mathfrak{R}\left\langle q^{\f}_{[c_{0}]},\mu\right\rangle 
		\end{align*} where $\mu\in T_{[c_{0}]}\T$ is the Beltrami differential denoting the derivative $\ddt [c_{t}]$.  }
	
	{\remark It is not known if the extremal length function is $C^{2}$ in general. See\cite{Royden1971,Liu2016,nate} for details on this regularity.That the Hubbard-Masur map may not be $C^{1}$ can be now seen from the above formula. 
	\subsection{Intersection of $q^{\f}$ and $q^{-\g}$ in $T^{*}\T$}\label{sectionintersectionofsection}\hfill\break\\
	Using the tools developed thus far we can state:	{\proposition \label{traninter}
		Let $(\mathsf{F},\mathsf{G})\in\mathcal{FMF}(\Sig)$ be a pair of  measured foliations that fill $S$ and $q^{\mathsf{F}},q^{-\mathsf{G}}: \mathcal{T}(S) \rightarrow T^{*}\mathcal{T}(S)$ be the associated sections defined before. Then their images in $T^{*}\T$ intersect uniquely and the projection of the intersection into $\T$ is the unique critical point of the function $\ext(\f)+\ext(\g):\T \rightarrow \mathbb{R}$. Moreover when $(\f,\g)\in \mathcal{FMF}_{0}(S)$ then the sections intersect transversely. }
	{\proof 
		Given, $(\mathsf{F},\mathsf{G}) \in \mathcal{FMF}(S)$, the sections $q^{\mathsf{F}}_{[c_{0}]},q^{-\mathsf{G}}_{[c_{0}]}$ intersect if and only if $q^{\f}_{[c_{0}]}=q^{-\g}_{[c_{0}]}=-q^{\g}_{[c_{0}]}$ from Remark \red{\ref{lelele}}. If $\mu$ be the Beltrami differential denoting $\ddt [c_{t}]$, then we have from Lemma \red{\ref{gardform}} that:
		
		\begin{align*}
			\ddt\ext_{[c_{t}]}(\mathsf{F})=  \mathfrak{R}(\langle q^{\mathsf{F}}_{[c_{0}]}, \mu \rangle)=\mathfrak{R}(\langle q^{-\mathsf{G}}_{[c_{0}]}, \mu \rangle),  \mu \in T_{[c]}\mathcal{T}(S).
		\end{align*}
		We can also consider $\g$ to be a measured foliation realised as the horizontal measured foliation and consider
		\begin{align*}
			\ddt\ext_{[c_{t}]}(\mathsf{G})= \mathfrak{R}(\langle q^{\mathsf{G}}_{[c_{0}]}, \mu \rangle),  \mu \in T_{[c]}\mathcal{T}(S).
		\end{align*}
		
		Hence $[c_{0}]$ is a critical point of the function $\ext(\mathsf{F})+\ext(\mathsf{G})$ if and only if $q^{\f}_{[c_{0}]}= -q^{\g}_{[c_{0}]}$. The existence and uniqueness of the critical point follows from Theorem \ref{GM} and Remark \ref{gola}. \\
		Now assume $(\f,\g)\in \mathcal{FMF}_{0}(S)$. If $dq^{\f},dq^{-\g}: T\T \rightarrow T T^{*}\T$ be the respective differentials then for transversality of intersection we need to show that if $dq^{\f}(\nu)=dq^{-\g}(\nu)$ for $\nu \in T\T$ then $\nu = 0$. Recalling the definitions of the sections $q^{\f},q^{-\g}$ this amounts to showing:
		{\lemma Consider a deformation of the type $(c_{t},q_{t}),t\geq 0$ with $(c_{0},q_{0})$ being the point of intersection of $q^{\f},q^{-\g}$ with $(\f,\g)\in \mathcal{FMF}_{0}(S)$. Let $\f_{t},\g_{t}$ be the horizontal and vertical measured foliations realised by $(c_{t},q_{t})$ and assume $\ddt \f_{t} = \ddt \g_{t}=0$. Then the deformation is trivial. }
		
		{\proof Recall that $q^{\f},q^{-\g}:\T \rightarrow T^{*}\T$ is a smooth map when $(\f,\g)$ are arational. For this, we consider $\pi:\widehat{S}\rightarrow S$ to be the canonical double cover branched over the zeroes of $q$  (see\red{\cite{Lanneau2003}}, Construction {$1.2$} and \S $5.4$ of\cite{Dumas2015}) such that $\pi^{*}(q)=\omega_{q}^{2}$ where $\omega_{q}$ is a holomorphic $1$ form on $(\widehat{S},[\widehat{c}])$ where $[\widehat{c}]=\pi^{*}[c]$.	 It follows from Lemma $2$ of\red{\cite{Lanneau2003}} that $([c],q)\mapsto ([\widehat{c}],\omega_{q}^{2})$ is a local embedding, so a deformation $([c_{t}],q_{t})$ in the generic stratum induces a deformation $([\widehat{c_{t}}],\omega_{q_{t}}^{2})\in Q(\widehat{\Sig})$ maintaining the same strata. Let $\widehat{\f_{t}}$ nd $\widehat{\g_{t}}$ be the horizontal and vertical foliations realised by $\omega_{q_{t}}^{2}$ on $(\widehat{\Sig},[\widehat{c_{t}}])$ with $\widehat{\f_{0}}=\widehat{\f}$ (resp. $\widehat{\g_{0}}=\widehat{\g}$) being the lift of $\f$ (resp. $\g$) in the double cover. Consider now  $\gamma$ to be a cycle in the  relative homology group $H^{-}_{1}(\widehat{\Sig},V_{\omega_{q}};\mathbb{C})$ where the latter is the eigenspace of $H_{1}(\widehat{\Sig},V_{\omega_{q}};\mathbb{C})$ consisting of cycles invariant under the involution of $\widehat{S}$ and the set $V_{\omega_{q}}$ denotes the set of zeroes of $\omega_{q}$. The real and imaginary part of the holonomy $\int_{\gamma}\omega_{q}$ are precisely the intersection number with the horizontal and vertical foliations of $\omega_{q}^{2}$. This gives us
			Period coordinates \begin{align*}
				\mathsf{per} &: Q_{0}(S)\hookrightarrow H^{1}_{-}(\widehat{S},V_{\omega_{q}};\mathbb{C})\cong \mathbb{R}^{12g-12}\\
				& q \mapsto (\gamma \rightarrow  \int_{\gamma} \omega_{q} )\end{align*} which is an immersion. Our assumption then translates to $\ddt i(\gamma,\widehat{\f_{t}})=0$ and $\ddt i(\gamma,\widehat{\g_{t}})=0$. This gives $\ddt \mathsf{per}(q_{t})(\gamma)=\ddt i(\gamma,\widehat{\g_{t}})+\textbf{i}i(\gamma,\widehat{\f_{t}})=0$, where we can assume $\gamma$ is fixed when one restricts to deformations maintaining strata. Since the period map is an immersion, it follows that this deformation is necessarily trivial. \qed \\\\ }
		Define now : {\definition For a pair $(\f,\g)$ that fill $\Sig$, we denote $\p(\f,\g)$ to be the critical point in $\T$ of the function $\ext(\f)+\ext(\g):\T \rightarrow \mathbb{R}$}.\\\\
		It is a simple observation from the definition that if the transverse measure of a foliation $\f$ is given by $|\mathfrak{R}\sqrt{q^{\f}}|$, then the corresponding holomorphic quadratic differential realising the measured foliation $t\f$ over the same Riemann surface structure on $\Sig$ is nothing but $t^{2}q^{\f}$ since then the transverse measure is given by $|\mathfrak{R}\sqrt{t^{2}q^{\f}}|$ which is equal to $t|\mathfrak{R}\sqrt{q^{\f}}|$. In the notation of the critical point $\p(\f,\g)$ this implies that:
		{\lemma\label{scalediff} If $(\f,\g)$ fill $S$, then $\p(t\f,t\g)=\p(\f,\g)$ where $t>0$.}
		{\proof Observe that $\ext(t\f)+\ext(t\g)= t^{2}(\ext(\f)+\ext(\g))$ from the definition of extremal length function and hence $\ext(t\f)+\ext(t\g)$ and $\ext(\f)+\ext(\g)$ have the same critical points. \qed} \\\\
		Also we have the observation that this point is uniquely determined by the second coordinate. That is:
		{\lemma\label{oneone} If $\mathsf{p}(\mathsf{F},\mathsf{G})=\mathsf{p}(\mathsf{F'},\mathsf{G}) $, then $\mathsf{F}= \mathsf{F'}$ in $\MF$.}
		{\proof Let $\p(\f,\g)=[c]$ be the unique point in $\T$ and $q\in Q(\Sig,c)$ be the unique holomorphic quadratic differential realising $(\f,\g)$ as its horizontal and vertical measured foliations respectively. For the pair $(\f',\g)$ we have that $\p(\f',\g)=p(\f,\g)=[c]$. Since on $Q(S,c)$ the choice of $q'$ realising $\g$ as its vertical measured foliation is unique from the theorem of Hubbard-Masur, we have that $q'=q$. But by definition $\f'$ is measure equivalent to $\mathsf{hor}_{[c]}(q')=\mathsf{hor}_{[c]}(q)=\f$.\qed}
		
		\subsection{Quotient of $Q(\Sig)$ under the action of $\mathbb{R}_{>0}$ and intersection of $[q^{\f}]$ and $[q^{-\g}]$ }\label{sectionquotientinter}\hfill\break \\
		There is a natural action of $\mathbb{R}_{>0}$ on $(Q(S, [c])-{0}) $ which sends every non-zero $q \in Q(S, [c]) $ to $ t^{2}q$,$ \forall t \in (0, \infty)$. We can thus define $Q^{1}(S,c)$ to be quotient $(Q(S, [c])-{0})/\mathbb{R}_{>0}$ under this action. Clearly $Q^{1}(S,c)$ is isomorphic to $ UT^{*}_{[c]}\mathcal{T}(S)$ from Proposition \ref{keylemma}, where the latter denotes the unit cotangent space at a point $[c] \in \T$. The next proposition is a similar result for the sections $[q^{\f}]$, which are the images of $q^{\mathsf{F}}$ under the quotient map. We can now address the main proposition of this section involving the intersection of the equivalence classes $[q^{\f}]$ and $[q^{-\g}]$ in $UT^{*}\T$ for a filling pair $(\f,\g)$: 
		{\proposition\label{line}
			Let $(\mathsf{F},\mathsf{G})\in \FMF$ be a pair filling measured foliations on $S$, then the projection of the intersection of the sections $[q^{\f}]$,$[q^{-\g}] $ in $ UT^{*}\mathcal{T}(S)$ onto $\mathcal{T}(S)$ is a  geodesic line for the Teichm{\"u}ller metric given by $t\mapsto \p(\sqrt{t}\f,\frac{1}{\sqrt{t}}\g)\in \T$ for $t>0$. Moreover, when $(\f,\g)\in \mathcal{FMF}_{0}(S)$ then the sections intersect transversely in $UT^{*}\T$.
		}
		
		{\proof
			
			We first note that if a pair $(\f,\g)$ fill $\Sig$ then so do the pairs $(t\f,\g)$,$(t\f,\frac{1}{t}\g)$ and $(\f,t\g)$ for any $t>0$. \\
			Let $[c_{t}]\in \T$ be an equivalence class of complex structures such that the two sections $q^{t\f},q^{-\g}$ meet over $[c_{t}]$. Then by definition we have $ t^{2}q^{\mathsf{F}}_{[c_{t}]}=q^{-\mathsf{G}}_{[c_{t}]}$ which is equivalent to $tq^{\f}_{[c_{t}]}=\frac{1}{t}q^{-\g}_{[c_{t}]}$. Since the foliation $t\f$ is realised by $t^{2}q^{\f}_{[c_{t}]}$ on the same complex structure, we have that $ q^{\sqrt{t}\mathsf{F}}_{[c_{t}]}=q^{-\frac{1}{\sqrt{t}}\mathsf{G}}_{[c_{t}]}$ for some $t>0$ at the point $[c_{t}]$. \\This is equivalent to the fact that $[c_{t}]$ is the unique critical point of the function $\ext(\sqrt{t}\mathsf{F})+\ext(\frac{1}{\sqrt{t}}\mathsf{G})$ since $(\sqrt{t}\mathsf{F},\frac{1}{\sqrt{t}}\mathsf{G})$ fill $S$. As $[c_{t}]$ is identified with $\mathsf{p}(\sqrt{t} \f, \frac{1}{\sqrt{t}}\g) $, the projection of the intersections is along 
			\begin{align*}
				\mathbb{R}_{>0} \rightarrow \T \\
				t \mapsto [c_{t}]=\p(\sqrt{t}\f,\frac{1}{\sqrt{t}}\g)
			\end{align*}
			So it now suffices to show that the path $t\rightarrow \p(t\f,\g)$ is a geodesic for the Teichm{\"u}ller metric $d_{\T}$ on $\T$. We first note that $[c]$ being the critical point $\ext(\f)+\ext(\g)$ for a filling pair $(\f,\g)$ also implies that $[c]$ is the critical point for the function
			\begin{align*}
				\ext(\f)\ext(\g): \T \rightarrow \mathbb{R}
			\end{align*} 
			where we use the fact that $\ext_{[c]}(\f)=\ext_{[c]}(\g)$ from corollary \ref{gorombhat}. So $\p(\sqrt{t}\f,\frac{1}{\sqrt{t}}\g)$ is a critical point for $\ext(\sqrt{t}\f)\ext(\frac{1}{\sqrt{t}}\g)$ and since $\ext_{[c]}(t\f)=t^{2}\ext_{[c]}(\f)$ we also have as a consequence that the point $\p(\sqrt{t}\f,\frac{1}{\sqrt{t}}\g)$ is a critical point for the function $\ext(\f)\ext(\g)$. Now it has been shown in\cite{Gardiner1991} that the set of critical points for the function $\ext(\f)\ext(\g)$ is a Teichm{\"u}ller geodesic line in $\T$ when $(\f,\g)$ fill $S$. Moreover, from Lemma \ref{oneone} the map $t\mapsto \p(t\f,\g)\in \T$ is injective. Finally, we observe that every critical point of $\ext(\f)\ext(\g)$ is also a critical point for $\ext(\alpha\f)+\ext(\beta \g)$ for some $\alpha, \beta >0$ and hence the the image of the map $t\rightarrow \p(t\f,\g)$ is the entire Teichm{\"u}ller geodesic. \\
			For transversality we can use Proposition \ref{traninter} as the pairs $q^{t\f}$ and $q^{-\g}$ intersect transversely, i.e, \begin{align*}
				T_{([c_{t}],q_{t})} T^{*}\T = T_{([c_{t}],q_{t})}(q^{t\f}(\T)) \bigoplus T_{([c_{t}],q_{t})} (q^{-\g} (\T))
			\end{align*} is true for all $t\geq 0$ and $(\f,\g)\in \mathcal{FMF}_{0}(S)$. The result follows when we take quotient.\qed   \\\\
			For a given pair $(\f,\g)\in \FMF$ we call  $t \mapsto \p(\sqrt{t}\f,\frac{1}{\sqrt{t}}\g)\in \T$ for $t>0$ as $\pp(\f,\g)$.

			\section{Necessary condition for paths with small filling measured foliations at infinity}\label{minsec}
			
			The goal of this section is to establish a necessary conditions that small differentiable  paths in $\QF$ starting from $\F$ should satisfy if the measured foliations at infinity are given by a filling pair $(t\fp,t\fm)\in \FMF$ at first order at $\F$. For this reason following\cite{Uhlenbeck1984}, we will study the curve $\beta_{([c],q)}(t^{2})\in \QF$, for $t>0$ small enough, which is parametrised by the data of the unique minimal surface it contains, i.e, the first fundamental form $I$ is in the conformal class $[c]\in \T$ and the second fundamental form $\II$ is given by $t^{2}\mathfrak{R}(q)$ for some $q\in T^{*}_{[c]}\T$. We will compute first-order estimates for Schwarzians at infinity for this path and determine that if the measured foliations at infinity for this path is indeed $(t\fp,t\fm)$ at first-order at $\F$ then $[c]$ is indeed the unique critical point for the functions $\ext(\fp)+\ext(\fm):\T\rightarrow \mathbb{R}$ and $q$ is the unique holomorphic differential we obtain from the theorem of Gardiner-Masur that realise $(\fp,\fm)$ on $[c]$.

			\subsection{Fundamental forms at infinity}\label{sectionfundadorms}\hfill\break\\
			Given a minimal surface in an almost-Fuchsian manifold $M$, we can consider the surfaces equidistant from it in $M$ at an oriented distance. These surfaces foliate the almost-Fuchsian manifold and we can then compute the associated first and second fundamental forms for these surfaces in terms of the data associated to the minimal embedding. We thus can formulate the following\red{\cite{Krasnov2008}}:
			{\lemma
				Let $S$ be a complete, oriented, smooth surface with principal curvatures in $(-1,1)$ immersed minimally into an almost-Fuchsian manifold homeomorphic to $S \times (-\infty,\infty)$ and let $(I,\II,B)$ be the associated data of the immersion. Then $\forall r \in \mathbb{R}$ the set of point $S_{r}$ at an oriented distance $r$ from $S$ is a smooth embedded surface with data $(I_{r},\II_{r},B_{r})$ where :
				\begin{enumerate}
					\item $I_{r}(x,y)=I((cosh(r)E+sinh(r)B)x,(cosh(r)E+sinh(r)B)y)$
					\item $\II_{r}=\frac{1}{2}\frac{dI_{r}}{dr}$
					\item $B_{r}=(cosh(r)E+sinh(r)B)^{-1}(sinh(r)E+cosh(r)B)$
				\end{enumerate}
				where $E$ is the identity operator and $S_{r}$ is identified to $S$ through the closest point projection.}\\
			
			The {fundamental forms at infinity} denoted as $I^{*},\II^{*}$ and introduced in\red{\cite{Krasnov2008}}, quantify the asymptotic behaviour of the quantities described above as $r \rightarrow \infty$. In particular, it estimates the data at the conformal class at infinity of an almost-Fuchsian manifold $M$ with respect to the, unique minimal surface with principal curvature in $(-1,1)$, it contains. \\\\
			Formally,
			\begin{align*}
				 I^{*}=\lim_{r \rightarrow \infty}2e^{-2r}I_{r} \hspace{1cm}\II^{*}=\lim_{r \rightarrow \infty}(I_{r}-\III_{r})
			\end{align*}. However the lemma above gives us explicit formulae to express the same in terms of $(I,\II,\III)$ and we use that to define:
			{\definition
				Adhering to the notations introduced above, the first fundamental form at infinity is given by the expression $I^{*}=\frac{1}{2}(I+2\II+\III)$ and 
				the second fundamental form at infinity is given by $\II^{*}=\frac{1}{2}(I-\III)$.
			}\\\\
			The pair $(I^{*},\II^{*})$ satisfy a modified version of Gauss equation at infinity (see\cite{Krasnov2007}) i.e, $\frac{1}{2}\trace(B^{*})=-K^{*}$ where $B^{*}$ is the shape operator associated to $I^{*}$ and $\II^{*}$. The Codazzi equation on the other hand, holds as it is by considering the Levi-Civita connection $\nabla^{*}$ compatible with $I^{*}$. The thing for importance to us is the expression for curvature associated to $I^{*}$ which we call $K^{*}$.\red{\cite{Krasnov2008}} further provide us with an expression for it using the data of the immersed minimal surface:
			{\lemma \label{formatinfinity}
				With the notation as above, 
				\begin{align*}
					K^{*}= \frac{K}{\deter(E+B)}=\frac{-1+\deter(B)}{1+\deter(B)}
				\end{align*} 
				where $K$ is the Gaussian curvature of the minimal immersion of $S$. }\\
			{\remark
				The second equality follows from the fact that $\trace(B)=0$ the immersion being minimal and $(I,\II)$ satisfy the Gauss-Codazzi equations.}\\
			
			In general $I^{*}$ need not be a hyperbolic metric. In fact,\red{\cite{Krasnov2008}} note that when multiplied by the correct conformal factor to take $I^{*}$ to the unique hyperbolic metric in its conformal class, the corresponding change in $\II^{*}$ is closely related to the Schwarzian derivative $\sigma$ associated to that end. So we have the following accounting for the change in $(\II^{*})_{0}$ when we apply a conformal change to $I^{*}$: 
			{\lemma\label{conformalchange} Let $I^{*}_{1}$ and $I^{*}_{2}$ be two metrics in the same conformal class at infinity such that $I^{*}_{2}=e^{2f}I^{*}_{1}$ for some smooth function $f$, then the traceless parts $(\II^{*}_{1})_{0}$ and $(\II^{*}_{2})_{0}$ are related as:
				\begin{align}\label{conformalchange}
					(\II^{*}_{2})_{0}-(\II^{*}_{1})_{0}=\hess_{I^{*}_{1}}(f)-df\otimes df + \frac{1}{2}||df||I^{*}_{1}-\frac{1}{2}(\Delta f)I^{*}_{1}
			\end{align}} In fact if we consider a holomorphic map $u:\Omega \rightarrow \mathbb{C}$ where $\Omega\subset \mathbb{C}$ then $\mathfrak{R}(\sigma(u))$ is precisely the term on the right hand side of the above equation when we consider $2f=\log(\frac{u'}{u})$. We thus have the following: (a geometric proof of which can also be found in Appendix $A$ of\red{\cite{Krasnov2008}})}:
		
		{\theorem\label{schlenker}
			If $I^{*}$ is hyperbolic, then $(\II^{*})_{0}=-\mathfrak{R}(\sigma)$, where $(\II^{*})_{0}$ denotes the traceless part of the second fundamental form at infinity and $\sigma$ is the Schwarzian at infinity.
		}\\\\
	We also note here that by the theorem above, upon uniformising $I^{*}$ to the hyperbolic metric $I^{*}_{h}$ in the conformal class and then computing $\II^{*}_{h}$ accordingly using Equation \ref{conformalchange}, we can see the map introduced in \ref{eq:schwarzian map} can be re-written as \begin{align*}
		\mathfrak{A}:	\QF \rightarrow T^{*}\mathcal{T}(\partial^{+}_{\infty}M)\\
		g\mapsto (I^{*}_{h},(\II^{*}_{h})_{0}).
	\end{align*}
		In the following sections we will use the parametrisation of almost-Fuchsian metric in terms of the data of $I$ and $\II$ of its unique minimal surface and compute $I^{*}$ and $\II^{*}$ at the two ends. For this, we will use the curve introduced by Uhlenbeck in\red{\cite{Uhlenbeck1984}} to prove that quasi Fuchsian metrics close enough to $\F$ admit a minimal surface with data given by a point $([c],sq)\in T^{*}\T$, for $s>0$ sufficiently small.\\

		\subsection{The curve $\beta_{([c],q)}(s)$ in $\QF$}\label{sectioncurveinQF}\hfill\break\\
		Let $([c],q)$ be a point in $ T^{*}\T$. As discussed in\red{\cite{Uhlenbeck1984}} we consider a smooth $1$-parameter curve $\beta_{([c],q)}(s)$ , $s \in [0,\epsilon)$  of almost-Fuchsian metrics starting from the Fuchsian locus which are given by the data 
		\begin{align*}
			\beta_{([c],q)}:\mathbb{R}_{>0}\rightarrow T^{*}\T\supset &\Omega\cong\AF\subset \QF\\
			s\mapsto ([c],s\mathfrak{R}(q))
		\end{align*}
		of the unique minimal surface such that $I$ is $e^{2u_{s}}h$ for some function $u_{s}:\Sig \rightarrow \mathbb{R}$, where $h$ denotes the unique hyperbolic metric in the conformal class $c$, and $\II=s\mathfrak{R}(q)$. At $s=0$, we have $u_{s}=0$ and $\beta_{([c],q)}(0)\in \F$. \\\\
		By Gauss equation, the pair $(e^{2u_{s}}h,s\mathfrak{R}(q))$ is the data of the minimal immersion if and only if $u_{s}$ is a solution for the following equation:
		\begin{align}\label{gausscodazzo}
			e^{-2u_{s}}(-\Delta_{h} u_{s}-1)= -1 + e^{-4u_{s}}s^{2}\deter_{h}(\mathfrak{R}(q)).
		\end{align}

		{\remark This is a reformulation of the Gauss equation $K_{g}=-1+\deter(B)$ for the pair $(e^{2u_{s}}h,s\mathfrak{R}(q))$. The left hand side comes from Lemma \red{\ref{formula}}. The right hand side comes by the formulae for change of basis for determinants.}\\\\
		It is then known from\cite{Uhlenbeck1984} (see also\cite{samuel,hodge} for this topic) that a unique solution exists for Equation \ref{gausscodazzo} which in terms of almost-Fuchsian metrics can be formulated as:
		
		{\proposition\label{citethisprop} 	For $0\leq s < \epsilon$, $\exists$ a unique almost-Fuchsian metric with a unique minimal surface whose $(I,\II)$ is given by the pair $([c],s\mathfrak{R}(q))$.} \\

		\subsection{First order estimations of measured foliations at infinity for the path $\beta_{([c],q)}(t^{2})$}\label{sectionfoliationsatinfinity}\hfill\break\\
		We will in fact do all the computations for the path $\beta_{([c],q)}(s)$ and perform a change of variable of $s$ to $t^{2}$ later on. This is done in order to account for the correct factor of the measured foliations at infinity at first order that we will compute eventually. Let us fix some notations: For a fixed $s>0$, the data of the minimal surface $S$ embedded into an almost-Fuchsian manifold $M$ can be expressed as 
		$I_{s}=e^{2u_{s}}h$ and $\II_{s}=s\mathfrak{R}(q)$.
		Since $B_{s}=I_{s}^{-1}\II_{s}$ and $\III_{s}(x,y)=I_{s}(B_{s}^{2}x,y)$, a simple computation in local orthonormal coordinates for $I_{s}$ show that $\III_{s}$ is equal to $-s^{2}e^{-2u_{s}}(\deter_{I_{s}}(\mathfrak{R}(q)))h$.
		Let the associated fundamental forms at infinity for this manifold be $I^{*}_{s},\II^{*}_{s}$ and the curvature at infinity be $K^{*}_{s}$. Further, let the Schwarzian at infinity associated to the two ends of $\beta_{([c],q)}(s)$ be called $\sigma^{s}_{+}$ and $\sigma^{s}_{-}$.\\\hspace*{.5 cm}
		Our goal first is to say that $I^{*}_{s}$ is hyperbolic at first order at $s=0$, so that we can apply a first order version of Theorem \red{\ref{schlenker}} relating the traceless part of $\II^{*}_{s}$ with the real part of $\sigma^{s}_{+}$.

		{\lemma\label{firstorderi}
			$I^{*}_{s}$ is hyperbolic at first order at $\F$ i.e the derivative of the curvature $K^{*}_{s}$ with respect to $s$ vanishes at $\F$ and $K^{*}_{0}=-1$ at $s=0$. Moreover, for this path $\ddtt  (\II^{*}_{s})_{0}=-\mathfrak{R}(q)$  .}
		{\proof
			First we note that at $s=0$ we are at the Fuchsian locus and from Lemma \ref{formatinfinity} we have $K^{*}_{0}=-1$.  Now observe that \begin{align*}
				\frac{d}{ds}\Bigr|_{\substack{s=0}} (\deter(B_{s})) = \frac{d}{ds}\Bigr|_{\substack{s=0}} s^{2}e^{-2u_{s}}\mathfrak{R}(q) = 0.
			\end{align*}Therefore using Lemma \ref{formatinfinity} \begin{align*}
				\frac{d}{ds}\Bigr|_{\substack{s=0}}	K^{*}_{s}= \frac{d}{ds}\Bigr|_{\substack{s=0}} \frac{-1+\deter(B_{s})}{1+\deter(B_{s})}=0.
			\end{align*}
			For the next part we first see that $u_{s}$ solves:
			
			\begin{align}\label{step1}
				e^{-2u_{s}}(-\Delta_{h} u_{s}-1)= -1 + e^{-4u_{s}}s^{2}\deter_{h}(\mathfrak{R}(q))
			\end{align}
			
			We define the non-linear map:
			\begin{align}
				F: W^{(2,2)}(S) \times [0,\infty)  \rightarrow   L^{2}(S) &&\\
				F(u_{s},s)=-\Delta_{h} u_{s} -1\label{step2} +e^{2u_{s}}-e^{-2u_{s}}s^{2}\deter_{h}(\mathfrak{R}(q)) 
			\end{align}
			where $W^{(2,2)}$ is the classical Sobolev space. The Fr{\'e}chet derivative is given by:
			\begin{align*}
				dF_{(u_{s},s)}(\dot u_{s},\dot s)=-\Delta \dot u_{s} +2\dot u_{s} e^{2u_{s}}+2\dot u_{s} e^{-2u_{s}}s^{2}\deter_{h}(\mathfrak{R}(q))-2e^{-2u_{s}}s \dot s \deter_{h}(\mathfrak{R}(q))
			\end{align*}It is clear that $u=u_{s}$ solves Equation (\ref{step1}) if and only if $(u_{s},s)$ is a solution for Equation (\ref{step2}). We now see the linearised operator with respect to $u$ of the function $F(u,s)$ which has the expression: 
			\begin{align*}
				L_{u_{s}}(\dot u_{s})=-\Delta_{h}\dot u_{s} +2\dot u_{s}(e^{2u_{s}}+e^{-2u_{s}}s^{2}\deter_{h}(\mathfrak{R}(q)).
			\end{align*}
			
			So, at the solution $(u_{s},s)=(0,0)$ we have that
			\begin{align*}
				L_{u_{0}}:W^{(2,2)}(S) \rightarrow L^{2}(S)\\ \dot u_{s} \mapsto -\Delta_{h} \dot u_{s}+ 2\dot u_{s}
			\end{align*} is a linear isomorphism of vector spaces see \cite{cmcdifian}, Lemma $3.4$. So  
			we can apply Implicit Function Theorem to get the solution curve $\gamma : [0,\epsilon)\rightarrow W^{(2,2)}(S)\times [0,+\infty)$ where $\gamma(s):=(u_{s},s)$  satisfies $F(u_s,s)=0,\forall s \in [0,\epsilon)$ ( see \cite{Uhlenbeck1984}). Now \begin{align*}
				dF_{(u_{s},s)}(\dot u_{s},\dot s)=L_{u_{s}}(\dot u_{s})-2e^{-2u_{s}}s\dot s \deter_{h}(\mathfrak{R}(q)).
			\end{align*}
			We have that $dF(\dot u_{s},\dot s)=0$ for this path so,
			\begin{align*} 
				\dot u_{s}&= L_{u_{s}}^{-1}(2 s \dot se^{-2u_{s}}\deter_{h}(\mathfrak{R}(q))\\
				& \implies \ddtt u_{s}=0.
			\end{align*}
			Now recall that $I^{*}_{s}=\frac{1}{2}(I_{s}+2\II_{s}+\III_{s})=\frac{1}{2}(e^{2u_{s}}h+2s\mathfrak{R}(q)-s^{2}e^{-2u_{s}}\deter_{I_{s}}(\mathfrak{R}(q))h)$. So taking derivative at $s=0$ gives us: 
			\begin{align}\label{japarishlekh}
				\frac{d}{ds}\Bigr|_{\substack{s=0}} I^{*}_{s}= \frac{1}{2}(2\dot u_{0}h+2\mathfrak{R}(q)+0)=\mathfrak{R}(q).
			\end{align}

			Now, by same computation observe that \begin{align*}
				\ddtt \II^{*}_{s}= \ddtt(I_{s}-\III_{s})=0.\end{align*} \\
			Now,\red{\cite{Krasnov2008}} further shows that that the mean curvature at infinity is expressed as $H^{*}_{s}=-K^{*}_{s}$. 
			Writing \begin{align*}
				&	\II^{*}_{s}=(\II^{*}_{s})_{0}+H^{*}_{s}I^{*}_{s}\\
				\implies \ddtt (-\II^{*}_{s})_{0} &= H^{*}_{0}\ddtt I^{*}_{s}+I^{*}_{0}\ddtt H^{*}_{s} \\
				\implies \ddtt(\II^{*}_{s})_{0}&= -\ddtt I^{*}_{s}.
			\end{align*} 
			From Equation \ref{japarishlekh} we have our claim.  \qed}\\ \hfill\break
		Now from Theorem \ref{schlenker} we know that if $I^{*}_{s}$ is hyperbolic then $(\II^{*}_{s})_{0}$ is equal to $-\mathfrak{R}(\sigma^{s}_{+})$ where $\sigma^{s}_{+}$ is the Schwarzian at the positive end at infinity. Moreover if we parametrise the quasi-Fuchsian space by the data of hyperbolic metric and Schwarzian at infinity at one end at infinity as in Equation \ref{eq:schwarzian map}, 
		then at the point $[c]\in \F$ of the Fuchsian locus we have a canonical decomposition of the tangent space $T_{[c]}(\QF)=T_{[c]}\T \bigoplus T^{*}_{[c]}\T $ where the first factor is the tangent to the Fuchsian locus denoting the derivative of the hyperbolic metric and the second factor is the derivative of the schwarzian at infinity at the Fuchsian locus. When considering the path $\beta_{([c],q)}(s)$ we have that $\ddtt \mathfrak{A}( \beta_{([c],q)}(s))=(\mathfrak{R}(q),-\mathfrak{R}(q))$. So, 
		
		{\lemma\label{firstorderschwarz}
			For the path $\beta_{([c],q)}(s)$, $\ddtt \sigma^{s}_{+}=q$ 
		}\\
		
		Note that we have done all the computation at one boundary component at infinity of $M$, which is almost-Fuchsian. However, recall that $M$ admits a foliation by surfaces "parallel" to the minimal surface, and the corresponding computation for the other component will differ by a sign. To be precise, $I^{*}_{s}=\frac{1}{2}(I_{s}-2\II_{s}+\III_{s})$ when we consider the component at the boundary at the other end at infinity (see \cite{Krasnov2007}). The rest of the computation follows as it is. Keeping this in mind we have:
		{\proposition\label{fluckit}
			For the path $\beta_{([c],q)}(s)$, $\ddtt \sigma^{s}_{\pm}=\pm q$  .}\\
		
		Upon a change of variable from $s$ to $t^{2}$, we will now show that the path $\beta_{([c],q)}(t^{2})$ is indeed a candidate for a path of almost-Fuchsian metrics with measured foliations at infinity given by the pair $(t\fp,t\fm)$ at first order. Denote the measured foliations at infinity for a metric in this path to be $\mathfrak{F}(\beta_{([c],q)})(t^{2})=(\f^{t}_{+},\f^{t}_{-})$. Here again a 1-parameter family of foliations $\f^{t}$ is said to be equivalent to a foliation $\f$ at first order, if for any given closed curve $\gamma$ on $\Sig$ 
		\begin{align*}
			\ddt i(\gamma,\f^{t} )=i(\gamma,\f) \\ 
			\implies i(\gamma,\f^{t})= ti(\gamma,\f)+o(t)
		\end{align*}
		Note by Proposition \red{\ref{fluckit}} $\sigma^{t^{2}}_{\pm}=\pm t^{2}q+o(t^{2})$ at first order at $t^{2}=0$ (or at $\F$) for the path $\beta_{([c],q)}(t^{2})$. So we need to show:
		
		{\lemma\label{firtorderfoliations} For any isotopy class of simple closed curve $\gamma$  on $S$ we have:
			\begin{align*}
				i(\gamma,\f^{t})-i(\gamma,\mathsf{hor}_{[c]}(t^{2}q))= o(t)
		\end{align*}	}
		{\proof 
			We just need to compute the following difference according to the definitions :
			\begin{align*}
				i(\gamma,\f^{t})-i(\gamma,\mathsf{hor}_{[c]}(t^{2}q))&= \inf_{\gamma}\int_{\gamma} | \mathfrak{I}\sqrt{\sigma_{t^{2}}}|- \inf_{\gamma}\int_{\gamma} |\mathfrak{I}\sqrt{t^{2}q} | \\
				&=\inf_{\gamma}\int_{\gamma} |\mathfrak{I} \sqrt{t^{2}q+o(t^{2})}|- \inf_{\gamma}\int_{\gamma} |\mathfrak{I}\sqrt{t^{2}q} | \\
				&=\inf_{\gamma}\int_{\gamma}|\mathfrak{I}\frac{\sigma_{t^{2}}-t^{2}q}{\sqrt{\sigma_{t^{2}}}+\sqrt{t^{2}q}}|
				=o(t).
			\end{align*}
			So, we have our claim. \qed \\
			\subsection{Necessary conditions for paths with given small filling measured foliations at infinity at first order}\label{sectioncondition}\hfill\break		
			So we see that for $0< t <\epsilon$ metrics in the path $\beta_{([c],q)}(t^{2})$ have that the measured foliations at infinity $(\f^{t}_{+},\f^{t}_{-})$ which at first order at the Fuchsian locus is given by the filling pair $(t\fp,t\fm)$. Secondly, notice that the point $([c],q)$ is the unique point associated to the filling pair $(\fp,\fm)$ via the Gardiner-Masur Theorem and $\beta_{([c],q)}(0)=[c]=\p(\fp,\fm)$ is the unique critical point for the functions $\ext(\fp)+\ext(\fm)$ by Proposition \red{\ref{traninter}}. These two points will precisely help us to formulate the condition we want paths with given first order behaviour of measured foliations at infinity to satisfy. 
			{\proposition \label{neccesary}  Let $(\fp,\fm)$ be a pair of measured foliations that fill $S$. Then there exists differentiable curve of quasi-Fuchsian metrics $t \mapsto \beta_{([c],q)}(t^{2})$, for $t\in[0,\epsilon)$, starting from the Fuchsian locus such that the image $\mathfrak{F}(\beta_{([c],q)}(t^{2}))\in \FMF$ is measure equivalent to $(t\fp,t\fm)$ at first order at $\F$. Moreover $[c]\in \T$ is the unique critical point of the function $\ext(\fp)+\ext(\fm): \T \rightarrow \mathbb{R}$ and $q\in T^{*}_{[c]}\T$ is the unique holomorphic quadratic differential realising $(\fp,\fm)$.}
			{\proof From  Proposition \ref{fluckit} and Lemma \ref{firtorderfoliations} we have that measured foliations at infinity for this path are given by the pair $(t\mathsf{hor}_{[c]}(q),t\mathsf{hor}_{[c]}(-q))$ at first order at $t=0$. The proposition is then a consequence of Proposition \ref{traninter}.\qed}

			\section{Uniqueness of paths with small filling foliations}\label{pathexists}	
			
			The goal of this section is to construct differentiable paths realising small pairs of measured foliations at infinity which are arational and filling, utilising the condition proved in Proposition \red{\ref{neccesary}} that they should satisfy. To do that first we will introduce the {blow-up} space $\widetilde{\QF}$ which we obtain by replacing $\F \subset \QF$ with its "unit normal bundle" $UN\mathcal{F}(\Sig)$. Following the strategy of\red{\cite{bonahon05}} 
			we then consider subsets of $\QF$ called $\Wfp$ (and $\Wfm$), defined as: {\definition For $\f\in \mathcal{MF}(S)$, define $\Wf^{+}\subset \QF$ (resp. $\Wf^{-}$) to be the set of quasi-Fuchsian metrics $g$ such that the foliation at the end at $+\infty$ (resp. $-\infty$) is $t\f$ for all $t\geq 0$.}\\\\Call $\widetilde{\mathcal{W}^{\pm}_{\f}}$ the image of $\Wf^{\pm}$ under the lift $\QF\rightarrow \widetilde{\QF}$. For $(\fp,\fm)\in \mathcal{FMF}_{0}(S)$ we will then show in the blow-up space $\widetilde{\Wfp}$ and $\widetilde{\Wfm}$ are submanifolds of $\widetilde{\QF}$ and that their boundaries $\partial \widetilde{\Wfp}$ and $\partial\widetilde{\Wfm}$ contained and intersecting in $\partial \widetilde{\QF}$ where $\partial \widetilde{\Wfp}\cap \partial \widetilde{\Wfm}$ intersect transversely and project onto the Teichm{\"u}ller geodesic line $\mathsf{P}(\fp,\fm)\in \T$ as defined in Proposition \red{\ref{line}}. 
			We then consider the map ${\pi}:{\Wfp}\cap{\Wfm}\rightarrow {\mathbb{R}^{2}}$, sending $g \in \widetilde{\Wfp}\cap\widetilde{\Wfm}$ to $(a,b)$ where $\mathfrak{F}(g)=(a\fp,b\fm)$ for some $a,b\geq 0$ by definition and lift the setting to the blow-up $\widetilde{\pi}: \widetilde{{\Wfp}}\cap \widetilde{\Wfm} \rightarrow \widetilde{\mathbb{R}^{2}}$ where the latter denotes the blow-up of $\mathbb{R}^{2}$ at the origin. The existence of paths with given small foliations then follows as we show that $\widetilde{\pi}$ at $\partial \widetilde{\Wfp} \cap \partial \widetilde{\Wfp}$ is a local diffeomorphism.

			\subsection{The normal bundle $N\F$ to $\F$ }\label{sectionnormal}\hfill\break\\
			First, let us recall that the Weil-Petersson metric endows $\T$ with a symplectic form $\omega_{WP}$ which is defined on the cotangent space as 
			\begin{align*}
				\omega_{WP}(.,.)= -\mathfrak{I}\left\langle .,.\right\rangle _{WP}.
			\end{align*} Moreover, $\T$ is endowed with an almost complex structure $J_{WP}$ such that $\left\langle q_{1},q_{2}\right\rangle_{WP} = \omega_{WP}(q_{1},J_{WP}(q_{2}))$ defined by $J_{WP}(q_{2})=\textbf{i}q_{2}$. Further, recall the notion of the {character variety} $\chi_{PSL_{2}(\mathbb{C})}$  which is an irreducible affine variety of complex dimension $6g-6$ and can be expressed as the GIT quotient:
			\begin{align*}
				\chi_{PSL_{2}(\mathbb{C})}:= \mathsf{Hom}(\pi_{1}(S),PSL_{2}(\mathbb{C}))//PSL_{2}(\mathbb{C})
			\end{align*} As each hyperbolic structure on $S$ is uniquely determined by the holonomy representation of $\pi_{1}(S)$ in to the group of orientation preserving isometries of $\mathbb{H}^{2}$, identified with $ PSL_{2}(\mathbb{R})$, the Fuchsian locus $\F(\cong \T)$ can be identified with a connected component of the set of {real points} in $\chi_{PSL_{2}(\mathbb{C})}$ (see\cite{Goldman1988}). Now, the group of orientation preserving isometries of $\mathbb{H}^{3}$ is identified with $PSL_{2}(\mathbb{C})$ and so the space $\QF$ is also identified with an open neighbourhood of $\F$ in $\chi_{PSL_{2}(\mathbb{C})}$ via discrete faithful representations from $\pi_{1}(S)\rightarrow PSL_{2}(\mathbb{C})$ which we can associate to a quasi-Fuchsian metric. This provides $\QF$ with a complex structure $J_{PSL_{2}(\mathbb{C})}^{2}=-1$ which also gives a decomposition of the tangent space at a point $[c]\in\F$ as $T_{[c]}\QF=T_{[c]}\F\bigoplus J_{PSL_{2}(\mathbb{C})}T_{[c]}\F
			$. This enables use to recall {Bers' Simultaneous Uniformisation Theorem}: 
			{\theorem[\red{\cite{bers1965moduli}}]\label{bersthm}
				The map $B:\QF \rightarrow \T \times \mathcal{T}(\overline{S})$ mapping a quasi-Fuchsian metric $ g \in \QF $ to the pair $B(g):=([c_{+}],[c_{-}])$ is biholomorphic with respect to the complex structure $J_{PSL_{2}(\mathbb{C})}$ of $\QF$ coming from the character variety and the complex structure $J_{WP}$ on $\T$.} \\\\
			It is also clear that $\F$ is the pre-image of the diagonal.  If $v\in T_{[c]}\QF$ is the tangent vector to the path $t\rightarrow g_{t}$ of quasi-Fuchsian metrics for $0\leq t <\epsilon$ at $t=0$ such that $g_{0}\in \F$ then the derivative of the Bers map at a point $[g_{0}]\in \T$  is given by  
			$
			d_{[g_{0}]} B(v) := (q_{1},q_{2})$
			where $q_{1},q_{2}$ are two holomorphic quadratic differentials in $Q(S,[g_{0}])$ denoting tangent vectors to $\T$ associated to the variation of the complex structures at two ends at infinity corresponding to the vector $v$. We thus have that $
			d_{[g_{0}]}B(J_{PSL_{2}(\mathbb{C})}v)=(J_{WP}(q_{1}),-J_{WP}(q_{2}))=(\textit{\textbf{i}}q_{1},-\textit{\textbf{i}}q_{2})
			$ where $\textit{\textbf{i}}^{2}=-1$ and the minus sign in the second factor is simply due to the opposite orientation of $S$.\\\hspace*{.5 cm}As $\F$ is identified with $\T$ by considering the unique hyperbolic metric $m$ in each conformal class $[c]$, when we consider a deformation of hyperbolic structures on $S$ the tangent vector $T_{m}\F\cong T_{m}\T$ is given by $\mathfrak{R}(q)$ for some $q\in Q(S,c)$ (see\cite{bookoftromba}). If one considers the variation of the hyperbolic metrics in the conformal classes associated to the two ends at infinity, then $d_{[c]}B(J_{PSL_{2}(\mathbb{C})}v)=(\mathfrak{R}(\textit{\textbf{i}}q_{1}),\mathfrak{R}(-\textit{\textbf{i}}q_{2}))\in T_{[c]}\F\times T_{[c]}\F$. We can thus define:
			\begin{definition}
				The {normal bundle} $N\F\rightarrow \F$ is the bundle whose fiber $N_{[c]}\F$ over each conformal class $[c]\in\T$ is the vector space isomorphic to the quotient $T_{[c]}\QF/T_{[c]}\F$.
			\end{definition}
			The fibers $N_{[c]}\F$ are $J_{PSL_{2}(\mathbb{C})}T_{[c]}\F$ where $J_{PSL_{2}(\mathbb{C})}$ is the almost complex structure of $\QF$. So $N_{[c]}\F$ is the set of tangent vectors $v_{([c],q)}$ such that $dB(J_{PSL_{2}(\mathbb{C})}v_{([c],q)})=(\mathfrak{R}(\textit{\textbf{i}}q), \mathfrak{R}(\textit{\textbf{i}}q))\in T\T \times T\T $ for $q\in Q(S,c)$.
			Now, let $v_{([c],q)}$ be the vector tangent to the path $\beta_{([c],q)}: [0,\epsilon) \rightarrow \QF$ at $\F$. That is to say:
			\begin{align*}
				v_{([c],q)}=\ddtt \beta_{([c],q)}(s)\in T_{[c]}\QF.
			\end{align*} So we formulate:
			{\proposition $v_{([c],q)}$ is an element of $N_{[c]}\F$. }
			{\proof As the first fundamental form at infinity $I^{*}_{s}$ remains hyperbolic at first order at $s=0$ it follows from Lemma \red{\ref{firstorderi}} that \begin{align*}
					dB(v_{([c],q)})=(\mathfrak{R}(q),-\mathfrak{R}(q)).
				\end{align*} Then, \begin{align*}
					dB(J_{PSL_{2}(\mathbb{C})}v_{([c],q)})=(\mathfrak{R}(\textit{\textbf{i}}q),\mathfrak{R}(\textit{\textbf{i}}q)).
				\end{align*}  So $J_{PSL_{2}(\mathbb{C})}v_{([c],q)}\in T_{[c]}\F$, i.e, tangent to the Fuchsian locus. The decomposition of the tangent space $T_{[c]}\QF= T_{[c]}\F \bigoplus J T_{[c]}\F $ at the Fuchsian locus then implies that $v_{([c],q)}\in N_{[c]}\F$.\qed\\
				
				\subsection{The blow-up $\widetilde{\QF}$ of $\QF$ at $\F$}\label{sectionblowup}\hfill\break\\
				For constructing the blow-up $\widetilde{\QF}$ consider again the bundle $N\F$ defined over $\F$. So we take the quotient of $N_{[c]}\F\setminus{0}$ by the action of $\mathbb{R}_{>0}$ called the unit normal bundle $UN_{[c]}\F$ and let $\overline{v_{([c],q)}}$ be the image of $v_{([c],q)}\in N_{[c]}\F$.  Consider now $\eta(N\F)\rightarrow UN\F$ to be the canonical differentiable line bundle and also we have a canonical linear map $\eta ({N\F)}\rightarrow \F $. We can show now that $\eta(N\F )\setminus {(\text 0-section)}\cong N\F\setminus \F$ is a diffeomorphism. Note that the zero section of $\eta(N\F)$ is again $UN\F$.\\\hspace*{.5 cm}
				Now let $\tau$ be a tubular map for $\F$ in $\QF$ and $\theta: \eta(N\F)\rightarrow N\F$ be the canonical map. The blow up $\widetilde{\QF}$ is the unique differentiable structure on $(\QF\setminus \F)\cup UN\F$ for which the inclusion map $\QF \setminus \F \subset \widetilde{\QF}$ and the map:
				\begin{align*}
					&\eta(N\F) \rightarrow \widetilde{\QF} \\
					&	v \mapsto \tau (\theta(v)) \hspace{.3cm}{\text{when}\hspace{.1cm} v\in \eta(N\F)\setminus{\text{(0-section)}}}\\
					&v \mapsto v \hspace{.3cm}\text{otherwise}
				\end{align*} are embeddings (see, for example\cite{Brocker1982}, pg. $128$). \\
				
				Moreover, it is a manifold with boundary $\partial \widetilde{\QF}$ which is $UN\F$. We observe that the natural inclusion $\QF\setminus\F\hookrightarrow \widetilde{\QF}$ lifts to $\QF\rightarrow \widetilde{\QF}$ by sending $[c]\in\F$ to $UN_{[c]}\F\in \partial \widetilde{\QF}$.\\\hspace*{.5 cm}
				Recall now the spaces of $\QF$, $\Wf^{+}$ and $\Wf^{-}$ for some $\f \in \mathcal{MF}(S)$.
				Since $\F \subset \Wf^{\pm}$, we have the natural inclusion $\Wf^{\pm}\setminus\F\hookrightarrow\QF\setminus\F$ which again lifts to a unique embedding $\Wf^{\pm} \rightarrow \widetilde{\QF}$ that we obtain by replacing a point $[c]\in \F \subset \Wf^{\pm} $ by the unique normal vector in $\overline{v_{([c],q)}}\in UN_{[c]}\F$ tangent to $\Wf^{\pm}$ at $[c]$. So by construction of the blow-up, for $t$ small enough $(\overline{v_{([c],q)}},t)\mapsto([c],t^{2}\mathfrak{R}(q))$ has $\f^{t}_{\pm}$ given by $t\f$ at first order at $\F$. By Lemma \red{\ref{firstorderschwarz}} and Theorem \red{\ref{hubbardmasur}},  $\overline{v_{([c],q^{\pm\f}_{[c]})}}$ is indeed that vector.
				We thus define: 
				{\definition $\widetilde{ \Wf^{+}}$ and $\widetilde{ \Wf^{-}}$ are the respective lifts of $\Wf^{+}$ and $\Wf^{-}$ into $\widetilde{\QF}$.}\\\\
				Having removed the Fuchsian locus which carry trivial Schwarzians from $\QF$, we will now parametrise elements in $\widetilde{\QF}$ by the data of the holomorphic quadratic differential being realised as the Schwarzian derivatives at the boundaries at infinity to show that $\widetilde{\Wf^{+}}$ and $\widetilde{\Wf^{-}}$ are submanifolds with boundary of $\widetilde{\QF}$. For this we first recall that the Schwarzians at infinity parametrise the $\mathbb{C}P^{1}$-structures on $\partial^{+}_{\infty}M$ and $\partial^{-}_{\infty}M$ (see \cite{Dumas}). More generally, if we denote the space of equivalence classes of $\mathbb{C}P^{1}$-structures on $\Sig$ under diffeomorphisms isotopic to the identity as $\mathcal{CP}(\Sig)$, then the Schwarzian derivative provides us parametrisation of the fibers of the forgetful map $\mathcal{CP}(S)\rightarrow \T$(\cite{Dumas}). So we formulate:
				{\lemma\label{graft} The Schwarzian parametrisation $\mathcal{S}:\QF \rightarrow T^{*}\mathcal{T}(\partial^{+}_{\infty}M)$ ( respectively for $\partial^{-}_{\infty}M$) introduced in Equation \ref{eq:schwarzian map}  is $C^{1}$.}   
				
				{\proof A quasi-Fuchsian metric on $M$ can be uniquely determined by the data of induced metric and measured bending lamination in the boundary of the convex core from\red{\cite{AFST_1996_6_5_2_233_0}} where we have that the map from $\QF \rightarrow\T \times \ML $ that associates the data of the unique pleated surface to the data of the quasi-Fuchsian metric is biholomorphic and so, smooth. Consider now the data $(m_{\pm},d_{m_{\pm}}(l(\lambda_{\pm})))$ which gives us a point in $T^{*}\T$, $d_{m_{\pm}}(l(\lambda_{\pm}))$ being the derivative of the length function for the measured lamination $\lambda_{\pm}$ computed at $m_{\pm}$. The claim is then a consequence of Theorem $4.1$ of\red{\cite{Dumas}}, originally due to Thurston and the main theorem in \cite{Krasnov2009}, which together state that the smooth \text{Grafting map} sending the data of the induced metric and measured bending lamination $(m_{\pm},\lambda_{\pm})\in \T \times \ML \cong T^{*}\T \rightarrow  \mathcal{CP}(S)$  on the boundary of the convex core $\partial^{\pm}_{\infty}\mathcal{CC}(M)$ to the data $([c_{\pm}],\sigma_{\pm})\in T^{*}\mathcal{T}(\partial^{\pm}_{\infty}M)$ at the boundary at infinity is a homeomorphism and $C^{1}$.\qed
					
					\subsection{Submanifolds $\widetilde{\Wfp}$, $\widetilde{\Wfm}$ and the intersection  $\partial \widetilde{ \Wfp}\cap\partial\widetilde{ \Wfm}$}\label{sectionblowupintersect}\hfill\break

					As Lemma \red{\ref{graft}} allows us to parametrise quasi-Fuchsian structures uniquely by the data of Schwarzian derivatives at the boundaries $\partial^{\pm}_{\infty}M$ we can thus proceed to discuss the following:
					
					{\proposition\label{5.2} For $\f \in \mathcal{MF}_{0}(S)$, the set ${\Wf^{+}}\setminus \F$ (resp. ${\Wf^{-}}\setminus  \F $) is a smooth submanifold of ${\QF}$ of dimension $\text{dim}(\T)+1$. In the blow-up $\widetilde{\QF}$, the lifts $\widetilde{\Wfp}$ and $\widetilde{\Wfm}$ are smooth submanifolds with the boundary $\partial \widetilde{\Wf^{+}}$ (resp. $\partial \widetilde{\Wf^{-}}$) contained in $\partial \widetilde{\QF}$.}
					{\proof We just treat the case of $\widetilde{ \Wf^{+}}$ as the same proof holds for $\widetilde{ \Wf^{-}}$ by symmetry. First we will show that $\Wf^{+}\setminus \F$ is a submanifold of $\QF\setminus\F$. Note that:
						\begin{align*}
							\QF \rightarrow \mathcal{CP}(\partial^{+}_{\infty}M)\hookrightarrow Q_{0}(\partial^{+}_{\infty}M)
						\end{align*}
						gives us an embedding of $\QF$ into an open subset of $Q_{0}(\partial^{+}_{\infty}M)$ by Lemma \ref{graft}.
						For a given $[c_{+}]$ we have an unique $\sigma^{\f}_{[c_{+}]}\in T^{*}\mathcal{T}(\partial^{+}_{\infty}M)$ realising $\f$ as the horizontal measured foliation at $\partial^{+}_{\infty}M$. This gives us the map $\sigma^{\f}:\mathcal{T}(\partial^{+}_{\infty}M) \rightarrow T^{*}\mathcal{T}(\partial^{+}_{\infty}M)$ which is identified with the map $q^{\f}$. Recall now that over the same complex structure $[c]$, $t\f$ is realised by $t^{2}\sigma^{\f}_{[c]}$. We see that $\Wf^{+}\setminus \F$ is locally embedded as $\mathbb{R}_{>0}\times \sigma^{\f}(\mathcal{T}(\partial^{\pm}_{\infty}M))$ in $\mathbb{R}^{12g-12}$ via the period coordinates of $\sigma^{\f}$. In other words, it is the image of the embedding:
						\begin{align*}
							\mathbb{R}_{> 0} \times \mathcal{T}(\partial^{+}_{\infty}M)\rightarrow Q_{0}(\partial^{+}_{\infty}M)\hookrightarrow \mathbb{R}^{12g-12}
						\end{align*} where the last inclusion is via the period coordinates associated to the dense stratum which gives us coordinate charts into $\mathbb{R}^{12g-12}$. Notice that the Fuchsian locus, corresponding to the zero section of $T^{*}\T$ has zero Schwarzian and thus the period associated is also zero. We also note that the smoothness of this submanifold is by virtue of the map $q^{\f}$ being real analytic when restricted to arational measured foliations. The dimension of the submanifolds being clearly $\text{dim}(\mathcal{T}(\partial^{\pm}_{\infty}M))+\text{dim}(\mathbb{R}_{>0})=\text{dim}(\T)+1$. \\
						Consider now the blow-up $\widetilde{Q_{0}(\partial^{+}_{\infty}M)}$ which is $(Q_{0}(S)\setminus\T)\cup Q^{1}_{0}(S)$ with the $C^{1}$ structure described in \S\S \ref{sectionblowup}. Since $(\overline{v_{[c],q^{\f}}})\in UN\F$ gets mapped to $\mathfrak{R}([q])\in UT\T$ under the isomorphism $UN\F\cong UT\T$ and $\mathfrak{R}(q)$ again corresponds to $[q]\in Q^{1}(S)$ by Weil-Petersson duality, we have an open embedding $\widetilde{\QF}\hookrightarrow \widetilde{Q_{0}(\partial^{+}_{\infty}M)}$. 
						Recall that $[c]\in \F\subset \Wf^{+}$ is associated to the unit normal vector $\overline{v_{([c],q^{\f}_{[c]})}}$ in $UN_{[c]}\F\subset \partial\widetilde{ \Wf^{+}}$ since $v_{([c],q^{\f}_{[c]})}$ is realised by the path $\beta_{([c],q^{\f}_{[c]})}(t^{2})$, and $\sigma^{t}_{+}$ for this path is indeed $t^{2}q^{\f}_{[c]}$ at first order at $\F$ by Lemma \red{\ref{firstorderschwarz}} and Lemma \red{\ref{firtorderfoliations}}. So $\partial\widetilde{\Wf^{+}}$ is contained in $Q^{1}_{0}(S)$, the boundary of $\widetilde{Q_{0}(S)}$.  \\
						For modifying the argument for $\widetilde{ \Wf^{-}}$ we need to consider the vector $q^{-\f}_{[c]}$ at $[c]\in \F\subset\Wf^{-}$, since the foliation at negative end at infinity for the path $\beta_{([c],q)}(t^{2})$ is given by $\mathsf{hor}_{[c]}(-q)$ at first order at $\F$ by Proposition \red{\ref{fluckit}} and Lemma \red{\ref{firtorderfoliations}}. The rest of the argument follows as it is and we have our claim. \qed}\\

					We can now claim the following:
					
					{\proposition\label{transbound} When $(\fp,\fm)\in \mathcal{FMF}_{0}(\Sig)$, $\partial\widetilde{\Wfp}$ and $\partial\widetilde{\Wfm}$ intersect transversely in $\partial\widetilde{\QF}$. Moreover, their intersection is equal to $(\mathsf{P}(\fp,\fm),\overline{v_{([c_{t}],q_{t})}})\ni UN\F $ where $[c_{t}]\in\T$ is the unique critical point of the function  $\ext(\sqrt{t}\fp)+\ext(\frac{1}{\sqrt{t}}\fm): \T \rightarrow \mathbb{R}$ and $q_{t}\in Q(S,c_{t})$ is the unique holomorphic quadratic differential realising them.}
					
					{\proof At $[c]\in \F$, $\partial\widetilde{\Wfp}$ is equal to $\overline{v_{([c],q^{\fp}_{[c]})}}\in UN_{[c]}\F$ and $\partial\widetilde{\Wfm}$ is given by the vector $\overline{v_{([c],q^{\fm}_{[c]})}}\in UN_{[c]}\F$. We know from \S \S \ref{sectionblowup} that $dB(J_{PSL_{2}(\mathbb{C})}\overline{v_{([c],q^{\fpm}_{[c]})}})\in UT_{[c]}\QF$ is given by $(\mathfrak{R}(\textit{\textbf{i}}[q^{\fpm}_{[c]}]),\mathfrak{R}(\textit{\textbf{i}}[q^{\fpm}_{[c]}]))\in UT_{[c]}\F\times UT_{[c]}\F$. So, $\partial\widetilde{\Wfp}$ intersects $\partial\widetilde{\Wfm}$ uniquely  at $[c]\in \F$ if and only if the sections $\mathfrak{R}([q^{\fp}])$ and $\mathfrak{R}([q^{\fm}])$ do in $UT\T$. Again, for some $[c]\in \T$, $\mathfrak{R}([q^{\fp}_{[c]}])=\mathfrak{R}([q^{\fm}_{[c]}])$ if and only if $[q^{\fp}_{[c]}]=[q^{\fm}_{[c]}]$ via the duality between $ T^{*}\T$ and $T\T$. So from Proposition \red{\ref{line}} we have that $[c]$ is the unique critical point $ \p(\sqrt{t}\fp,\frac{1}{\sqrt{t}}\fm)$ of the function $\ext(\sqrt{t}\fp)+\ext(\frac{1}{\sqrt{t}}\fm):\T\rightarrow \mathbb{R}$ for some $t>0$. From proposition \ref{line} we see that the boundaries $\partial \widetilde{ \Wfp}$ and $\partial\widetilde{ \Wfm}$ intersect along the image of the Teichm{\"u}ller geodesic line $t \mapsto \p(\sqrt{t}\fp,\frac{1}{\sqrt{t}}\fm)$ under the section map $\F\rightarrow UN\F$. The transversality of their intersection follows from that of the submanifolds $[q^{\fp}]$ and $[q^{-\fm}]$ shown in Lemma \ref{line}.
						\qed}}\\\\
				Define now the map 
				\begin{align*}
					\pi : \Wfp \cap \Wfm \rightarrow \mathbb{R}\times\mathbb{R}
				\end{align*}
				which sends $g\in \Wfp \cap \Wfm$ to the pair $(a,b)$ such that
				$\mathfrak{F}(g)=(a\fp,b\fm)$ by definition of $\mathfrak{F}$. Observe that under this map, $\F$ gets mapped to $\left\lbrace 0\right\rbrace :=(0,0)$ and for $(\fp,\fm)\in \mathcal{FMF}_{0}(S)$ the map $\pi$ is smooth. So, if $\widetilde{\mathbb{R}^{2}}$ be the blow-up of $\mathbb{R}^{2}$ at the origin, then $\pi$ lifts to a smooth map 
				\begin{align*}
					\widetilde{\pi}:\widetilde{\Wfp}\cap\widetilde{\Wfm}\rightarrow \widetilde{\mathbb{R}^{2}}.
				\end{align*}
				Here $\widetilde{\mathbb{R}^{2}}$ is the set $\mathbb{R}^{2}\setminus\left\lbrace 0\right\rbrace \cup UT_{\left\lbrace 0\right\rbrace }\mathbb{R}^{2}$ where $UT_{\left\lbrace 0\right\rbrace }\mathbb{R}^{2}$ is the quotient of the tangent space at origin under the action of $\mathbb{R}_{>0}$. 
				So we have 
				{\proposition \label{localdiffeo} 
					For a pair $(\fp,\fm)\in \mathcal{FMF}_{0}(S)$, the map $\widetilde{\pi}$ is a local diffeomorphism near $\partial\widetilde{\Wfp}\cap \partial\widetilde{\Wfm}$ onto its image.}
				
				{\proof 
					We want to show that the map $\widetilde{\pi}$ has a solution at the Fuchsian locus, is invertible at that point and subsequently apply implicit function theorem. For this we show that $\widetilde{\pi}$ is a local immersion and local submersion at $p\in\partial\widetilde{ \Wfp}\cap\partial\widetilde{\Wfm}$, i.e to prove that \begin{align*}
						d_{p}\widetilde{\pi}: T_{p}\widetilde{\Wfp} \cap T_{p}\widetilde{\Wfm} \rightarrow T_{\widetilde{\pi}(p)}\widetilde{\mathbb{R}^{2}}
					\end{align*} is injective and surjective. Note that when restricting to $(\fp,\fm)$ in arational pairs, this map is indeed smooth as the submanifolds $\widetilde{\mathcal{W}}^{\pm}_{\f_{\pm}}$ are for $\f\in \mathcal{MF}_{0}(S)$.\\
					If for some $v \in  T_{p}\widetilde{\Wfp} \cap T_{p}\widetilde{\Wfm} $ we have that $d_{p}\widetilde{\pi}(v)=\left\lbrace 0\right\rbrace $ then we want to show first that $v$ is in the intersection of the tangent spaces to the boundary $T_{p}\partial \widetilde{\Wfp}\cap T_{p}\partial \widetilde{\Wfm}$. Let $m_{+}:\Wfp\rightarrow [0,\infty)$ be the map such that 
					$m_{+}(g)= t$ for some $g\in \Wfp$ which has measured foliation at the boundary at infinity given by $t\fp$. This induces a map $\widetilde{m_{+}}: \widetilde{\Wfp}\rightarrow [0,\infty)$ in the blow-up space as well. Observe that if we analogously define a map $\widetilde{m_{-}}:\widetilde{\Wfm} \rightarrow [0,\infty)$ then $\widetilde{\pi}:=(\widetilde{m_{+}},\widetilde{m_{-}})$. So if $v\in \mathsf{Ker}(d\widetilde{\pi})$ then $v\in \mathsf{Ker}(d\widetilde{m_{+}})\cap\mathsf{Ker}(d\widetilde{m_{-}})$ then this implies that $v\in T_{p}\partial \widetilde{\Wfp} \cap T_{p} \partial\widetilde{\Wfm}$.
					
					\begin{figure}
						\centering
						\includegraphics[width=0.7\linewidth]{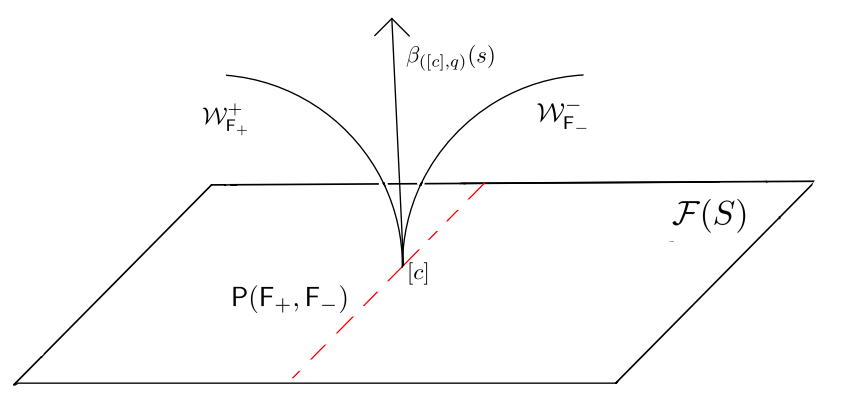}
						\caption{Schematic diagram of the path $\beta_{([c],q)}(s)$ leaving the Fuchsian locus $\F$ from the point $[c]$ along the direction of the normal vector $v_{([c],q)}$ with $\Wfp$ and $\Wfm$ intersecting at $\F\subset \QF$ prior to blow-up procedure for an arational pair $(\fp,\fm)$ which fills $S$ and the dashed line representing the Teichm{\"u}ller geodesic $\pp(\fp,\fm)$.}
						\label{fig:pathintersect}
					\end{figure}
					
					Now let $\kappa: \widetilde{\mathbb{R}^{2}_{\geq 0}}\rightarrow \mathbb{R}_{>0}\times\mathbb{R}_{\geq 0}  $ be the chart which sends $(x,y) \mapsto (\frac{x}{y},y)$. Recall from Proposition \red{\ref{transbound}} that $ \partial\widetilde{\Wfp}\cap \partial\widetilde{\Wfm} $ has been shown to be the lift of the line $\mathsf{P}(\fp,\fm)$ in $\T\cong\F$ by the section $\F\rightarrow UN\F$ sending $\mathsf{P}(\fp,\fm) \ni [c]\mapsto ([c],\overline{v_{([c],q^{\fp}_{[c]})}})=([c],\overline{v_{([c],q^{-\fm}_{[c]})}}) \in UN_{[c]}\F$.\\ For a fixed $t>0$, let $[c]$ be the critical point of $\ext(\sqrt{t}\fp)+\ext(\frac{1}{\sqrt{t}}\fm)$ which is equivalent to being the critical point of $\ext(t\fp)+\ext(\fm)$. Let $s\mapsto \widetilde{g_{s}} \in \widetilde{\Wfp}\cap \widetilde{\Wfm}$ be a differentiable path such that $\widetilde{g_{0}}$ is the  $([c],\overline{v_{([c],q^{\fp}_{[c]})}})=([c],\overline{v_{([c],q^{-\fm}_{[c]})}})\in \partial\widetilde{ \Wfp}\cap \partial \widetilde{ \Wfm}$ and suppose $\widetilde{g_{s}}$ in turn descends to a curve $g_{s} \in \QF$ with $g_{0}=[c]$ under the projection $\widetilde{\QF} \rightarrow \QF$. As $\widetilde{g_{s}}\in\widetilde{\Wfp} \cap \widetilde{\Wfm}$ we have that $\pi(\mathfrak{F}(\widetilde{g_{s}}))=(\widetilde{a}(s)\fp,\widetilde{b}(s)\fm)$ which again descend to two smooth functions $1$-parameter functions $a(s),b(s)\in \mathbb{R}_{\geq 0}$ such that $\mathfrak{F}(g_{s})=(a(s)\fp,b(s)\fm)$ and $a(0)=b(0)=0$. Moreover by definition $g_{s}$ is normal to $\F$ at $g_{0}$ and along the direction ${v_{([c],q^{\fp}_{[c]})}}={v_{([c],q^{-\fm}_{[c]})}}$.\\
					This brings us back to the case of Proposition \red{\ref{neccesary}} where we have a path starting from $\F$, normal to $\F$ and with specified first order behavior of the measured foliations at infinity given by a pair that fills $S$. Thus $g_{s}$ is a path of the type $\beta_{([c],q)}(t^{2})$ where $g_{0}=[c]=\mathsf{p}(a'(0)\fp,b'(0)\fm)$, the critical point for the function $\ext(a'(0)\fp)+\ext(b'(0)\fm)$. Again by assumption $[c]=\mathsf{p}(t\fp,\fm)$; so $\frac{a'(0)}{b'(0)}=t$, as $\mathsf{p}(\fp,\fm) $ is unique for a filling pair $(\fp,\fm)$ up to  scaling by $t$ (see Remark \red{\ref{scalediff}}). \\
					So we see that \begin{align*}
						\kappa \circ \widetilde{\pi}([c],\overline{{v_{([c],q^{\fp}_{[c]})}}})=\lim_{s \rightarrow 0}\kappa \circ \widetilde{\pi}(\tilde{g_{s}})=\lim_{s \rightarrow 0} \kappa \circ \pi(g_{s})=\lim_{s \rightarrow 0} \kappa \circ (a(s),b(s))=\lim_{s \rightarrow 0}(\frac{a(s)}{b(s)},b(s))=(t,0)
					\end{align*} \\
					This shows that if $v \in  T_{p}\widetilde{\Wfp} \cap T_{p}\widetilde{\Wfm}$ with $d_{p}\widetilde{\pi}(v)=0$ then $v$ is zero. Hence $d_{p}\widetilde{\pi}$ is injective at $\partial\widetilde{ \Wfp}\cap \partial\widetilde{ \Wfm}$.\\
					So the map $\widetilde{\pi}$ is a local immersion into $\widetilde{\mathbb{R}^{2}}$ at the points $p \in \partial\widetilde{\Wfp}\cap \partial\widetilde{\Wfm} $. Also $d_{p}\widetilde{\pi}$ is surjective at ${\partial\widetilde{\Wfp}\cap \partial\widetilde{\Wfm} }$ because the domain is a $2$ dimensional real manifold being the boundary of $\widetilde{ \Wf^{+}}\cap\widetilde{ \Wf^{-}}$ and the image is the boundary of the $2$ dimensional real manifold $\widetilde{\mathbb{R}^{2}}$ being $UT_{\left\lbrace 0\right\rbrace }\mathbb{R}^{2}$. \\
					So we proved that $\widetilde{\pi}$ is a local diffemorphism in a neighbourhood of $\partial \widetilde{ \Wfp}\cap\partial{\widetilde \Wfm}$. \qed   }
				
				We can now address the main proposition of this section which proves Theorem \ref{thm1.1}: 
				{\proposition\label{pathexist}
					Let $(\fp,\fm)$ be a pair of arational measured foliations that fill $S$ and let $\mathsf{p}(\fp,\fm)\in \F$ be the critical point of the function $\ext(\fp)+\ext(\fm)$. Then for $ t\in [0,\epsilon) $ there exists a unique smooth curve $t \mapsto g_{t}\in \QF$, with $g_{0}=\mathsf{p}(\fp,\fm)$, such that the $\mathfrak{F}(g_{t})=(t\fp,t\fm)$ for all $t\in[0,\epsilon)$.   
				}
				{\proof By the preceding Proposition there exists a $ t\rightarrow\widetilde{g_{t}} \in \widetilde{\Wfp}\cap\widetilde{\Wfm}$ smooth curve in an open neighbourhood of $\partial\widetilde{ \Wfp}\cap\partial\widetilde{ \Wfm}$ such $b \circ \widetilde{\pi} \circ \mathfrak{F}(\widetilde{g_{t}})=(t,t)$ for $t\in (0,\epsilon)$ with $b$ being the blow-up map $b:\widetilde{\mathbb{R}^{2}} \rightarrow \mathbb{R}^{2}$ mapping $UT_{\left\lbrace 0\right\rbrace }(\mathbb{R}^{2})$ to the origin and identity on the rest. The result then follows as $\widetilde{g_{t}}$ descends to $g_{t} \in \QF$ with $\mathfrak{F}(g_{t})=(t\fp,t\fm)$ for $t\in[0,\epsilon)$. \qed   }
				.
					\section{Interpretation in half-pipe geometry}\label{halfpipe}
				\subsection{Quasi-Fuchsian half-pipe $3$ manifolds}
				We will now give an interpretation of our result in quasi-Fuchsian {Half pipe} $3$ manifolds that we describe following\red{\cite{Danciger2013}}. To describe the space $\mathbb{HP}^{3}$, we will switch our viewpoint to the projective model for $\mathbb{H}^{3}$ in this section. Consider $\mathbb{RP}^{3}\subset \mathbb{R}^{4}$ with the group $PGL_{4}(\mathbb{R})$ being its isometry group.
				Consider now $\mathbb{H}^{3}$ as a subset of $\mathbb{RP}^{3}$. To be precise, consider $\mathbb{R}^{4}$ with the diagonal form given by the matrix 
				\begin{align*}
					\eta_{t}=\begin{bmatrix}
						-1&  &  & \\ 
						&  1&  & \\ 
						&  &  1& \\ 
						&  &  & t^{2}
					\end{bmatrix}
				\end{align*}	
				where $t \geq 0$. Each form $\eta_{t}$ define a convex region in $\mathbb{X}_{t}\subset \mathbb{RP}^{3}$ given by the relation 
				\begin{align*}
					x^{T}\eta_{t}x= -x_{1}^{2}+x_{2}^{2}+x_{3}^{2}+t^{2}x_{4}^{2}<0
				\end{align*}
				For each $t,$ $\mathbb{X}_{t}$ is a homogeneous subspace of $\mathbb{RP}^{3}$ which is preserved by the group $G_{t}$ of linear transformations that preserve $\eta_{t}$. With these notations, $\mathbb{H}^{3}=\mathbb{X}_{+1}$ and $G_{+1}=PO(3,1)\cong PSL_{2}(\mathbb{C})$. \\\\
				Moreover, define $\g_{t}: \mathbb{X}_{+1}\rightarrow \mathbb{X}_{t}$ as 
				\begin{align*}
					\begin{bmatrix}
						1&  &  & \\ 
						&  1&  & \\ 
						&  &  1& \\ 
						&  &  & t^{-1}
					\end{bmatrix}
				\end{align*}
				and this gives an isomorphism between $\mathbb{X}_{+1}$ and $ \mathbb{X}_{t}$ . Moreover $\g_{t}$ conjugates $PO(3,1)$ to $G_{t}$. 
				Notice further the co-dimension $1$ space $\mathbb{P}^{3}$ defined by $x_{4}=0$ and $-x_{1}^{2}+x_{2}^{2}+x_{3}^{2}<0$ is a totally geodesic copy of $\mathbb{H}^{2}$ and is contained in $\mathbb{X}_{t}$ for all $t$ and $g_{t}$ fixes $\mathbb{P}^{3}$ pointwise. \\\hspace*{.5 cm}
				For $t>0$ we now consider a $1$-parameter family of quasi-Fuchsian structures on $M\cong \Sig \times \mathbb{R}$. So, we have a family of developing maps and holonomy representations given by:
				\begin{align*}
					\mathcal{D}_{t}: \widetilde{\Sig} \rightarrow \mathbb{H}^{3}\cong X_{+1} \\
					\rho_{t}: \pi_{1}(\Sig)\rightarrow PO(3,1) \cong PSL_{2}(\mathbb{C})
				\end{align*}  
				Assume further that for $t=0$, $\mathcal{D}_{0}$ gives us a submersion of $\widetilde{\Sig}$ onto $\mathbb{P}^{3}=\mathbb{H}^{2}$. That is the coordinate $x_{4}$ converges to a zero function. $\rho_{t}$ then converges to $\rho_{0}$ whose image lies in the subgroup $PO(2,1)\cong SL(2,\mathbb{R}). $ \\\hspace*{.5 cm}
				Apply now the rescaling map to obtain the developing map $g_{t}\mathcal{D}_{t}:\widetilde{\Sig}\rightarrow\mathbb{X}_{t}$, so that the holonomy representation is given by $g_{t}\rho_{t}g_{t}^{-1}$. Suppose that $t \rightarrow 0$ then $g_{t}\mathcal{D}_{t}$ converges to a local diffeomorphism $\mathcal{D}:\widetilde{\Sig}\rightarrow \mathbb{X}_{0}$ and if $\rho_{\mathcal{D}}:\pi_{1}(\Sig)\rightarrow PGL_{4}(\mathbb{R})$ is the limit of the holonomy $\rho_{t}$ as $t \rightarrow 0$ then $\mathcal{D}$ is equivariant with respect to $\rho_{\mathcal{D}}$ . To be precise, for $\gamma \in \pi_{1}(\Sig)$ if $\rho_{t}$ is of the form: 
				
				\begin{align}\rho_{t}=\begin{pmatrix}
						A(t) &w(t) \\ 
						v(t)& a(t)
				\end{pmatrix}\end{align}
				where $A \in PO(2,1)\cong PSL(2,\mathbb{R})$ and $w(t),v(t)^{T} \in \mathbb{R}^{3} $, then we have
				\begin{align}\label{halfpipething}
					\lim_{t \to 0}g_{t}\rho_{t}(\gamma)g_{t}^{-1}=\lim_{t \to 0} \begin{pmatrix}
						A(t) &tw(t) \\ 
						v(t)/t& a(t)
					\end{pmatrix}=\begin{pmatrix}
						A(0) &0 \\ 
						v'(0)& 1
					\end{pmatrix}=\rho_{\mathcal{D}}
				\end{align}
				So we have the following: 
				
				{\definition\label{quasifuchsianhalfpipe} A half-pipe structure on $\Sig\times \mathbb{R}$ is a $(G_{HP^{3}},\mathbb{HP}^{3})$ structure where $\mathbb{HP}^{3}=\mathbb{X}_{0}$ and $G_{HP}$ is the subgroup of $PGL_{4}(\mathbb{R})$ of matrices with the form
					$\begin{pmatrix}
						A &0 \\ 
						v& \pm 1
					\end{pmatrix}$
					where $A \in O(2,1)$ and $v^{T}\in \mathbb{R}^{3}$.}\\
				We also define:
				{\definition\label{mcqhp} Any path $\rho_{t}$ of representations into $PSL_{2}(\mathbb{C})$ satisfying Equation (\ref{halfpipething}) is said to be compatible at first order at $t=0$ with $\rho_{\mathcal{D}}$}.
				\\ As observed in\red{\cite{Danciger2013}}, $G_{HP^{3}}\cong \mathbb{R}^{2,1} \rtimes O(2,1) $, where an element of the form $\begin{pmatrix}
					A &0 \\ 
					v& \pm 1
				\end{pmatrix}$ can be interpreted as an infinitesimal deformations of the the hyperbolic structure given by $A\in PO(2,1)$ and along the direction $v$ normal to $PO(2,1)$ into $PO(3,1)$. Passing onto quotients, we see that quasi-Fuchsian half-pipe $3$-manifolds are precisely obtained by infinitesimal deformations in $\QF$ starting from the point $[c]\in \F$ along a direction $v_{([c],q)} \in N_{[c]}\mathcal{F}(\Sig)$. So we define:
				{\definition\label{hpqf} $M^{HP}_{c,q}$ is the half-pipe quasi-Fuchsian structures whose holonomy representation into $PGL_{4}(\mathbb{R})$ is compatible at first order at $t=0$ with the holonomy $\rho_{t}$ associated to quasi-Fuchsian metrics in $\beta_{([c],q)}(t)$ in the sense of Definition \red{\ref{mcqhp}}.}
				\subsection{Half-pipe Schwarzians and their measured  foliations}\label{hpschwarzandfoli}\hfill\break\\
				Recall again that $v_{([c],q)}$ is the tangent vector to the path $\beta_{([c],q)}(t)\in \QF$ associated to which we have unique minimal immersions of $\Sig$ for each $t<\epsilon$ with immersion data $I_{t}\in [c]$ and $\II_{t}=t\mathfrak{R}(q)$. So for each $t$ we have a $\mathcal{D}_{t}$ and $\rho_{t}$ in the sense above and a half-pipe structure as the limit when $t$ goes to $0$. There is also an analogous notion for half-pipe for the second fundamental form and shape operator in half-pipe geometry that follows from \cite{Fillastre2019}. So we want to study the limit of these immersion data as $t\rightarrow 0$ and use the following lemma: 
				{\lemma [\red{\cite{Fillastre2019}}]Let $\sigma_{t}$ be a $C^{2}$ family of minimal immersions of $\mathbb{H}^{2}$ into $\mathbb{H}^{3}$, such that $\sigma_{0}$ is an embedding of $\mathbb{H}^{2}$. Let
					$\sigma = \lim_{t \to 0}G_{t}\circ\sigma_{t}$
					be the rescaled immersion in $\mathbb{HP}^{3}$. Then:
					\begin{itemize}
						\item  The first fundamental form of $\sigma$ coincides with the first fundamental form of $\sigma_{0}$:
						\begin{align*}
							I(v,w)=\lim_{t \to 0}I_{t}(v,w)
						\end{align*}
						\item The second fundamental form of $\sigma$ is the first derivative of the second fundamental form of $\sigma_{t}$:
						\begin{align*}
							\II(v,w)=\lim_{t \to 0}\frac{\II_{t}(v,w)}{t}
						\end{align*}
						\item The shape operator $B$ of $\sigma$ is the first derivative of the shape operator $B_{t}$ of $\sigma_{t}$:
						\begin{align*}
							B(v)= \lim_{t \to 0}\frac{B_{t}}{t}
						\end{align*}
				\end{itemize} }\hfill\break
				We immediately have the following for $M^{HP}_{c,q}$. 
				
				{\proposition\label{8,7}  The half-pipe manifold $M^{HP}_{c,q}$ contains a smooth minimal surface with immersion data uniquely given by $I\in[c]$ and $\II=\mathfrak{R}(q)$.  }
				{\proof Since $M^{HP}_{c,q}$ appears as a limit of the quasi-Fuchsian structures defines by the path $\beta_{([c],q)}(t)$ we apply the lemma directly. The induced metric of the minimal immersion of $\Sig$ into $M^{HP}_{c,q}$ is given by $I= \lim_{t \to 0}I_{t}\in [c]$. So we see that $\II=\lim_{t \to 0}\frac{\II_{t}}{t}=\mathfrak{R}(q)$.}\\\\
				So we will introduce an analogous notion for Schwarzian at infinity for half-pipe manifolds that is quite natural with the tools we have developed so far and with our definition of $M^{HP}_{cq}$.
				
				{\definition\label{definitionofhpsch}  The positive (resp. negative) half pipe Schwarzian at infinity associated to $M^{HP}_{cq}$ is defined as the derivative at $\F$ of the Schwarzian derivatives at the positive (resp. negative) end at infinity for quasi-Fuchsian metrics in the path $\beta_{([c],q)}(t)$ for $t$ small enough.}\\\\\
				From Lemma \red{\ref{firstorderschwarz}} So we have 
				{\proposition The positive and negative half-pipe Schwarzians at infinity for $M^{HP}_{c,q}$ are $q$ and $-q$.  }\\\\
				We can now again consider the horizontal measured foliation $\pm \f$ associated to $\pm q$ on $[c]$ and obtain our Theorem \ref{thm1.3} by an application of Theorem \red{\ref{GM}}:
				
				{\theorem\label{hpschth} Any pair $(\fp,\fm)\in \FMF$ can be uniquely realised as the horizontal foliations of the positive and negative half pipe Schwarzians at infinity associated to quasi-Fuchsian half-pipe manifold. Moreover, $M^{HP}_{cq}$ defined before is the unique one realising $(\fp,\fm)$, where $([c],q)\in T^{*}\T$ is the unique point realising $(\fp,\fm)$. }\\\\
				This can be seen as a {first-order interpretation of Theorem of Gardiner-Masur}  as half-pipe quasi-Fuchsian manifolds correspond to points in $T^{*}\T$ by Proposition \red{\ref{8,7}} via the minimal surface they contain.

				\bibliographystyle{abbrv}
				\bibliography{main}

			\end{document}